\documentclass[final]{siamltex}

\usepackage{amssymb}
\usepackage{color}
\usepackage{epsfig}
\usepackage{graphicx}
\usepackage{makeidx}
\usepackage{multicol}
\usepackage{enumerate}
\usepackage{amsfonts}
\usepackage{amsmath}
\usepackage{amsbsy}
\usepackage{algorithm}
\usepackage{algpseudocode}

\newcommand{\oprank}[1]{\mbox{{rank}$_{\otimes}$($#1$)}}

\newcommand{\mlrank}[1]{\mbox{{rank}$_{\boxplus}$($#1$)}}

\newdimen\pHeight
\newdimen\pLower
\newdimen\pLineWidth
\newdimen\pKern
\newdimen\pIR
\newsavebox{\Cbox}
\newsavebox{\vertCmplx}
\newdimen\Cheight
\newdimen\Cwidth
\sbox{\Cbox}{\rm C}
\sbox{\vertCmplx}{\rule[\pLower]{\pLineWidth}{\Cheight}}
\sbox{\Cbox}{\usebox{\Cbox}\kern\pKern\usebox{\vertCmplx}}
\wd\Cbox=\Cwidth

\def\R{{\rm I\kern\pIR R}}

\mathcode `:="003A
\newcommand{\ins}[1]{\mbox{${#1} \in \R    $}} 

\newcommand{\inv}[2]{\mbox{${#1} \in
\R^{\hspace*{-.01in} #2}$}}  

\newcommand{\inm}[3]
   {\mbox{${#1} \in \R^{\hspace*{-.015in}{#2} \times
{#3}}\hspace*{-.05in}$ }}

\newcommand{\inT}[5]
   {\mbox{${\cal #1} \in \R^{\hspace*{-.015in}{#2} \times
{#3} \times {#4} \times {#5}}\hspace*{-.05in}$ }} 

\newcommand{\inTens}[4]
   {\mbox{${\cal #1}_{#2} \in \R^{\hspace*{-.015in}{#3} \times \cdots \times {#4}} \hspace*{-.05in}$ }}

\newcommand{\norm}[1]{\mbox{$\|\: #1 \: \|$}}

\newcommand{\RR} { \mathbb{R} }   

\newcommand{\myskip}   { \medskip }

\newcommand{\T}{\raisebox{1pt}{$ \:\otimes \:$}}

\mathcode `:="003A

\title{Block Tensors and Symmetric Embeddings}
\author{Stefan Ragnarsson\thanks{Center for Applied Mathematics, Cornell University 
        Ithaca, NY 14853, {\tt str23@cornell.edu}}
\and Charles F. Van Loan\thanks{Department of Computer Science, Cornell University,
        Ithaca, NY 14853, {\tt cv@cs.cornell.edu}. Both authors
are supported in part by NSF contract DMS-1016284.}
}
\begin{document}
\maketitle
\begin{abstract}
Well known connections exist between the singular value decomposition of
a matrix $A$ and the Schur decomposition of its symmetric embedding
$\mbox{\bf sym}(A) = ([\:0\: A\:;\: A^{T}\:0\:])$. In particular, if $\sigma$ is a singular value of
$A$ then $+\sigma$ and $-\sigma$ are eigenvalues of the symmetric embedding. The top and
bottom halves of $\mbox{\bf sym}(A)$'s eigenvectors are singular vectors for $A$. 
Power methods applied to $A$ can be related to power methods applied to $\mbox{\bf sym}(A)$. The
rank of $\mbox{\bf sym}(A)$ is twice the rank of $A$.
In this paper we develop similar connections for tensors by building on L-H. Lim's variational approach to tensor singular values and vectors. We show how to embed a general order-$d$ tensor $\cal A$ into
an order-$d$  symmetric tensor $\mbox{\bf sym}({\cal A})$. Through the embedding we
relate  power methods for $\cal A$'s singular
values to power methods for $\mbox{\bf sym}({\cal A})$'s eigenvalues.
Finally, we connect the multilinear and outer product rank
of $\cal A$ to the multilinear and outer product rank of $\mbox{\bf sym}({\cal A})$. 
\end{abstract}


\pagestyle{myheadings}
\thispagestyle{plain}
\markboth{S. RAGNARSSON and C.F. VAN LOAN}{Symmetric Embedding of Tensors}

\section{Introduction}

If \inm{A}{n_{1}}{n_{2}}, then there are well-known connections between its singular value decomposition (SVD) and
the eigenvalue and eigenvector properties  of the symmetric matrix
\begin{equation}
\inm{\mathbf{sym}(A) \;=\; \left[ \begin{array}{cc} 0 & A \\ A^{T} & 0 \end{array}\right]}{(n_{1}+n_{2})}{(n_{1}+n_{2})}.
\end{equation}
If $A = U\Sigma V^{T}$ is the SVD of $A$, then for $k=1:\mbox{rank}(A)$
\begin{equation}
\left[ \begin{array}{cc} 0 & A \\ A^{T} & 0 \end{array}\right]
\left[ \begin{array}{c} u_{k} \\ \pm v_{k} \end{array} \right]
\;=\;
\pm \sigma_{k} \left[ \begin{array}{c} u_{k} \\ \pm v_{k} \end{array} \right]
\end{equation}
where $u_{k} = U(:,k)$, $v_{k} = V(:,k)$, and  $\sigma_{k} = \Sigma(k,k)$.
Another way to connect $A$ and $\mathbf{sym}(A)$ is through the Rayleigh quotients
\begin{equation}
\phi_{A}(u,v) \;= \; \frac{u^{T}Av}{\norm{u}_{2}\norm{v}_{2}} \;=\;
\left(
\sum_{i_{1}=1}^{n_{1}}\sum_{i_{2}=1}^{n_{2}} A(i_{1},i_{2})u(i_{1})v(i_{2})
\right) \Bigg/ \left( \norm{u}_{2}\norm{v}_{2}\right)
\end{equation}
 and
\begin{equation}
\rule{.16in}{0in} \phi_{A}^{(sym)}(x) \;=\; \frac{1}{2} \frac{x^{T}Cx}{x^{T}x} \;=\;\frac{1}{2}
\left( \sum_{i_{1}=1}^{N} \sum_{i_{2}=1}^{N} C(i_{1},i_{2})x(i_{1})x(i_{2})\right) \Bigg/
\norm{x}_{2}^{2}
\end{equation}
where  \inv{u}{n_{1}}, \inv{v}{n_{2}}, $N = n_{1}+n_{2}$, \inv{x}{N}, and $C = \mathbf{sym}(A)$.
If $x$ is a stationary vector for $\phi_{A}^{(sym)}$, then $u = x(1:n_{1})$ and $v = x(n_{1}+1:n_{1}+n_{2})$
render a stationary value for $\phi_{A}$. See \cite[p.448]{gvl}.

In this paper we discuss these notions as they apply to tensors.
An order-$d$ tensor $\inTens{A}{}{n_{1}}{n_{d}}$
is a real $d$-dimensional array $ {\cal A}(1:n_{1},\ldots,1:n_{d})$
where the index range in the $k$-th mode is from 1 to $n_{k}$. 
The idea of embedding a general tensor into a larger symmetric tensor having the same order is developed in \S2.  This requires having a facility with block tensors. Fundamental  orderings, unfoldings, and multilinear summations
are discussed in \S3 and used in \S4 where we characterize  various  multilinear Rayleigh quotients and their stationary values and vectors.  This builds on the variational approach to tensor singular values developed in \cite{lim2}. In \S5 we provide a symmetric embedding  analysis of several higher-order power methods for tensors that have recently been proposed \cite{regalia,regalia2,lath3,lath4,kolda4}. Results that 
relate the multilinear and outer product ranks of a tensor to the corresponding ranks of its symmetric embedding are presented in \S6. A brief conclusion section follows.

Before we proceed with the rest of the paper, we use the case of third-order tensors
to preview some of the main ideas and to establish notation. (The busy reader already familiar with basic tensor computations and notation may safely skip to \S2.)
The starting point is to define the trilinear Rayleigh quotient
\begin{equation}
\phi_{\cal A}(u,v,w) \;=\; 
\left(
\sum_{i_{1}=1}^{n_{1}}\sum_{i_{2}=1}^{n_{2}}\sum_{i_{3}=1}^{n_{3}}
{\cal A}(i_{1},i_{2},i_{3})
u(i_{1})v(i_{2})w(i_{3})
\right) \Bigg/
\left( \norm{u}_{2}\:\norm{v}_{2} \:\norm{w}_{2}\right)
\end{equation}
where  \inv{\cal A}{n_{1}\times n_{2} \times n_{3}},\inv{u}{n_{1}}, \inv{v}{n_{2}}, and \inv{w}{n_{3}}. Calligraphic characters
are used for tensors: ${\cal A}(i_{1},i_{2},i_{3})$ is entry $(i_{1},i_{2},i_{3})$ of $\cal A$.

The singular values and vectors of
${\cal A}$ are the critical values and vectors of $\phi_{\cal A}$ as formulated in \cite{lim2}. 
A simple expression for the gradient $\nabla \phi_{\cal A}$ is made possible by
unfolding ${\cal A} = (a_{ijk})$ in each of its three modes and aggregating the $u$, $v$, and $w$ vectors with the Kronecker product. To illustrate, suppose $n_{1} = 4$, $n_{2} = 3$, and $n_{3} = 2$ and define the
{\em modal unfoldings} ${\cal A}_{(1)}$, ${\cal A}_{(2)}$, and ${\cal A}_{(3)}$ by 

\begin{eqnarray}
{\cal A}_{(1)} &=& \left[ 
\begin{array}{cccccc}
a_{111} \!&\! a_{121} \!&\! a_{131} \!&\! a_{112} \!&\!  a_{122} \!&\! a_{132} \\
a_{211} \!&\! a_{221} \!&\! a_{231} \!&\! a_{212} \!&\!  a_{222} \!&\! a_{232} \\
a_{311} \!&\! a_{321} \!&\! a_{331} \!&\! a_{312} \!&\!  a_{322} \!&\! a_{332} \\
a_{411} \!&\! a_{421} \!&\! a_{431} \!&\! a_{412} \!&\!  a_{422} \!&\! a_{432} 
\end{array}
\right] \nonumber\\
& &  \nonumber \\
{\cal A}_{(2)} &=&\left[ 
\begin{array}{cccccccc}
a_{111} \!&\! a_{211} \!&\! a_{311} \!&\! a_{411} \!&\! a_{112} \!&\! a_{212} \!&\! a_{312} \!&\! a_{412} \\
a_{121} \!&\! a_{221} \!&\! a_{321} \!&\! a_{421} \!&\! a_{122} \!&\! a_{222} \!&\! a_{322} \!&\! a_{422} \\
a_{131} \!&\! a_{231} \!&\! a_{331} \!&\! a_{431} \!&\! a_{132} \!&\! a_{232} \!&\! a_{332} \!&\! a_{432} 
\end{array}
\right]\\
& &  \nonumber \\
{\cal A}_{(3)} &=&\left[ 
\begin{array}{cccccccccccc}
a_{111} \!&\! a_{211} \!&\! a_{311} \!&\! a_{411} \!&\! a_{121} \!&\! a_{221} \!&\! a_{321} \!&\! a_{421} \!&\! a_{131} \!&\! 
a_{231} \!&\! a_{331} \!&\! a_{431} \\
a_{112} \!&\! a_{212} \!&\! a_{312} \!&\! a_{412} \!&\! a_{122} \!&\! a_{222} \!&\! a_{322} \!&\! a_{422} \!&\! a_{132} \!&\! 
a_{232} \!&\! a_{332} \!&\! a_{432} 
\end{array}\right]. \nonumber
\end{eqnarray}
The columns of these matrices are {\em fibers}. A fiber of a tensor is obtained by fixing
all but one of the indices. For example, the third column of the unfolding
\[
{\cal A}_{(1)} \;=\;
\left[ \begin{array}{cccccc}
{\cal A}(:,1,1)  & {\cal A}(:,2,1)  & {\cal A}(:,3,1)  & {\cal A}(:,1,2)  & {\cal A}(:,2,2)  & {\cal A}(:,3,2)  
\end{array}
\right]
\]
is the fiber
\[
{\cal A}(:,3,1) \;=\; \left[ \begin{array}{c} {\cal A}(1,3,1) \\ {\cal A}(2,3,1) \\ {\cal A}(3,3,1) \\ {\cal A}(4,3,1)
\end{array} \right]
\]
obtained by fixing the 2-mode index at 3 and the 3-mode index at 1. It is necessary to
specify the order in which
the fibers appear in a modal unfolding. The choice exhibited in (1.6) has the property that
\begin{equation}
\sum_{i_{1}=1}^{n_{1}}\sum_{i_{2}=1}^{n_{2}}\sum_{i_{3}=1}^{n_{3}}
{\cal A}(i_{1},i_{2},i_{3})
u(i_{1})v(i_{2})w(i_{3})
\;=\;
\left\{
\begin{array}{l}
u^{T} \,{\cal A}_{(1)} \: w \T v \\
v^{T} \,{\cal A}_{(2)} \: w \T u  \rule{0pt}{14pt}\\
w^{T} {\cal A}_{(3)} \: v \T u \rule{0pt}{14pt}
\end{array}
\right.
\end{equation}
which makes it easy to specify  the stationary vectors of $\phi_{\cal A}$.
If $u$, $v$, and $w$ are unit vectors, then the gradient of $\phi_{\cal A}$ is given by
\begin{equation}
\nabla \phi_{\cal A}(u,v,w) \;=\;
\left[ \begin{array}{c}
{\cal A}_{(1)}\: w \T v \\
{\cal A}_{(2)}\: w \T u \rule{0pt}{14pt}\\
{\cal A}_{(3)}\: v \T u \rule{0pt}{14pt}
\end{array}
\right]
\;-\;
\phi_{\cal A}(u,v,w) \left[ \begin{array}{c} u \\ v \rule{0pt}{14pt}\\
\rule{0pt}{14pt} w
\end{array}\right].
\end{equation}
We remark that if $\cal A$ is an order-2 tensor, then (1.8) collapses to the
familiar matrix-SVD equations
$Av = \sigma u$ and $A^{T}u = \sigma v$.

A central contribution of this paper revolves around the tensor version of the {\bf sym} matrix
(1.1) and the associated Rayleigh quotient  $\phi_{\cal A}^{(sym)}$ that is defined in (1.4). Just as {\bf sym}-of-a-matrix sets up a symmetric block matrix whose entries are either zero or matrix transpositions,
{\bf sym}-of-a-tensor sets up a symmetric block tensor whose entries are either zero or  a tensor transposition.

If \inv{\cal A}{n_{1}\times n_{2} \times n_{3}}, then there are $6 = 3!$
possible transpositions identified by the notation ${\cal A}^{<\,[i\:j\:k]\,>}$
where $[i\:j\:k]$ is a permutation of $[1\:2\:3]$:
\begin{equation}
{\cal B} \;=\; \left\{ \begin{array}{l} {\cal A}^{<\:[1\;2\;3]\:>} \\ {\cal A}^{<\:[1\;3\;2]\:>} \rule{0pt}{14pt} \\
{\cal A}^{<\:[2\;1\;3]\:>}\rule{0pt}{14pt} \\
 {\cal A}^{<\:[2\;3\;1]\:>}\rule{0pt}{14pt} \\
  {\cal A}^{<\:[3\;1\;2]\:>} \rule{0pt}{14pt}\\ 
  {\cal A}^{<\:[3\;2\;1]\:>} \rule{0pt}{14pt} \end{array} \right\}
\qquad
\Longrightarrow \qquad
\left\{ \begin{array}{l}
{\cal B}(i,j,k) \\ 
{\cal B}(i,k,j) \rule{0pt}{14pt} \\
{\cal B}(j,i,k) \rule{0pt}{14pt}  \\ 
 {\cal B}(j,k,i) \rule{0pt}{14pt} \\
  {\cal B}(k,i,j) \rule{0pt}{14pt} \\ 
  {\cal B}(k,j,i) \rule{0pt}{14pt} \end{array} \right\}
\;=\; {\cal A}(i,j,k) 
\end{equation}
for $i=1:n_{1},\;j=1:n_{2},\; k=1:n_{3}.$

The symmetric embedding of a 3rd-order tensor results in a 3-by-3-by-3 block tensor,
a kind of  Rubik's cube built from 27 (possibly non-cubical) boxes.
If \inv{\cal A}{n_{1}\times n_{2} \times n_{3}} and $N = n_{1}+n_{2} + n_{3}$, then
\inv{\mbox{\bf sym}(\mathcal{A})\;=\;C}{N \times N \times N} is the 3-by-3-by-3 block tensor
whose $ijk$ block is specified by
\begin{equation} 
{\cal C}_{[i\:j\:k]} \;=\; 
\left\{
\begin{array}{ll} {\cal A}^{<\,[i\:j\:k]\,>} & \mbox{if $[i\:j\:k]$ is a permutation of
$[1\:2\:3]$}\\
&  \\ 
\inv{0}{n_{i}\times n_{j} \times n_{k}}
& \mbox{otherwise}.
\end{array}
\right.
\end{equation}
See {\sc Fig 1.1}. The blocks in a block tensor such as $\cal C$ can be specified using the colon notation.
For example, if $n_{1} = 4$, $n_{2} = 3$ and $n_{3} = 2$, then

\begin{equation}
\begin{array}{lclcl} 
{\cal C}_{[1\:2\:3]} &=& {\cal C}(1:4,5:7,8:9) & = & 
\inv{{\cal A}^{<\:[1\;2\;3]\:>}\;}{n_{1}\times n_{2} \times n_{3}} \\
{\cal C}_{[1\:3\:2]} &=& {\cal C}(1:4,8:9,5:7)& = & 
\inv{{\cal A}^{<\:[1\;3\;2]\:>}\;}{n_{1}\times n_{3} \times n_{2}} \rule{0pt}{14pt}\\
{\cal C}_{[2\:1\:3]} &=& {\cal C}(5:7,1:4,8:9)& = & 
\inv{{\cal A}^{<\:[2\;1\;3]\:>}\;}{n_{2}\times n_{1} \times n_{3}}\rule{0pt}{14pt}\\
{\cal C}_{[2\:3\:1]} &=& {\cal C}(5:7,8:9,1:4)& = & 
\inv{{\cal A}^{<\:[2\;3\;1]\:>}\;}{n_{2}\times n_{3} \times n_{1}} \rule{0pt}{14pt}\\
{\cal C}_{[3\:1\:2]} &=& {\cal C}(8:9,1:4,5:7)& = & 
\inv{{\cal A}^{<\:[3\;1\;2]\:>}\;}{n_{3}\times n_{1} \times n_{2}} \rule{0pt}{14pt}\\
{\cal C}_{[3\:2\:1]} &=& {\cal C}(8:9,5:7,1:4)& = & 
\inv{{\cal A}^{<\:[3\;2\;1]\:>}\;}{n_{3}\times n_{2} \times n_{1}} \rule{0pt}{14pt}
\end{array}\; .
\end{equation}
We will prove in section 2.3 that the tensor ${\cal C}$ is in fact symmetric.

\begin{figure}
\scalebox{.75}{
\begin{picture}(300,300)(-10,0)
\put(90,170){\makebox(0,0){\LARGE ${\cal C}(:,:,1)$}}
\put(10,10){\thicklines \line(1,0){120}}
\put(10,50){\thicklines \line(1,0){120}}
\put(10,90){\thicklines \line(1,0){120}}
\put(10,130){\thicklines \line(1,0){120}}
\put(10,10){\thicklines \line(0,1){120}}
\put(50,10){\thicklines \line(0,1){120}}
\put(90,10){\thicklines \line(0,1){120}}
\put(130,10){\thicklines \line(0,1){120}}
\put(10,130){\thicklines \line(1,1){20}}
\put(50,130){\thicklines \line(1,1){20}}
\put(90,130){\thicklines \line(1,1){20}}
\put(130,130){\thicklines \line(1,1){20}}
\put(30,150){\thicklines \line(1,0){120}}
\put(130,130){\thicklines \line(1,1){20}}
\put(130,90){\thicklines \line(1,1){20}}
\put(130,50){\thicklines \line(1,1){20}}
\put(130,10){\thicklines \line(1,1){20}}
\put(150,30){\thicklines \line(0,1){120}}
\put(70,30){\makebox(0,0){\large ${\cal A}^{<[3 2 1]>}$}}
\put(110,70){\makebox(0,0){\large ${\cal A}^{<[2 3 1]>}$}}

\put(250,210){\makebox(0,0){\LARGE ${\cal C}(:,:,2)$}}
\put(170,50){\thicklines \line(1,0){120}}
\put(170,90){\thicklines \line(1,0){120}}
\put(170,130){\thicklines \line(1,0){120}}
\put(170,170){\thicklines \line(1,0){120}}
\put(170,50){\thicklines \line(0,1){120}}
\put(210,50){\thicklines \line(0,1){120}}
\put(250,50){\thicklines \line(0,1){120}}
\put(290,50){\thicklines \line(0,1){120}}
\put(170,170){\thicklines \line(1,1){20}}
\put(210,170){\thicklines \line(1,1){20}}
\put(250,170){\thicklines \line(1,1){20}}
\put(290,170){\thicklines \line(1,1){20}}
\put(190,190){\thicklines \line(1,0){120}}
\put(290,170){\thicklines \line(1,1){20}}
\put(290,130){\thicklines \line(1,1){20}}
\put(290,90){\thicklines \line(1,1){20}}
\put(290,50){\thicklines \line(1,1){20}}
\put(310,70){\thicklines \line(0,1){120}}
\put(190,70){\makebox(0,0){\large ${\cal A}^{<[3 1 2]>}$}}
\put(270,150){\makebox(0,0){\large ${\cal A}^{<[1 3 2]>}$}}

\put(410,250){\makebox(0,0){\LARGE ${\cal C}(:,:,3)$}}
\put(330,90){\thicklines \line(1,0){120}}
\put(330,130){\thicklines \line(1,0){120}}
\put(330,170){\thicklines \line(1,0){120}}
\put(330,210){\thicklines \line(1,0){120}}
\put(330,90){\thicklines \line(0,1){120}}
\put(370,90){\thicklines \line(0,1){120}}
\put(410,90){\thicklines \line(0,1){120}}
\put(450,90){\thicklines \line(0,1){120}}
\put(330,210){\thicklines \line(1,1){20}}
\put(370,210){\thicklines \line(1,1){20}}
\put(410,210){\thicklines \line(1,1){20}}
\put(450,210){\thicklines \line(1,1){20}}
\put(350,230){\thicklines \line(1,0){120}}
\put(450,210){\thicklines \line(1,1){20}}
\put(450,170){\thicklines \line(1,1){20}}
\put(450,130){\thicklines \line(1,1){20}}
\put(450,90){\thicklines \line(1,1){20}}
\put(470,110){\thicklines \line(0,1){120}}
\put(390,190){\makebox(0,0){\large ${\cal A}^{<[1 2 3]>}$}}
\put(350,150){\makebox(0,0){\large ${\cal A}^{<[2 1 3]>}$}}

\end{picture}}
\caption{The Symmetric Embedding of an Order-3 Tensor}
\end{figure}
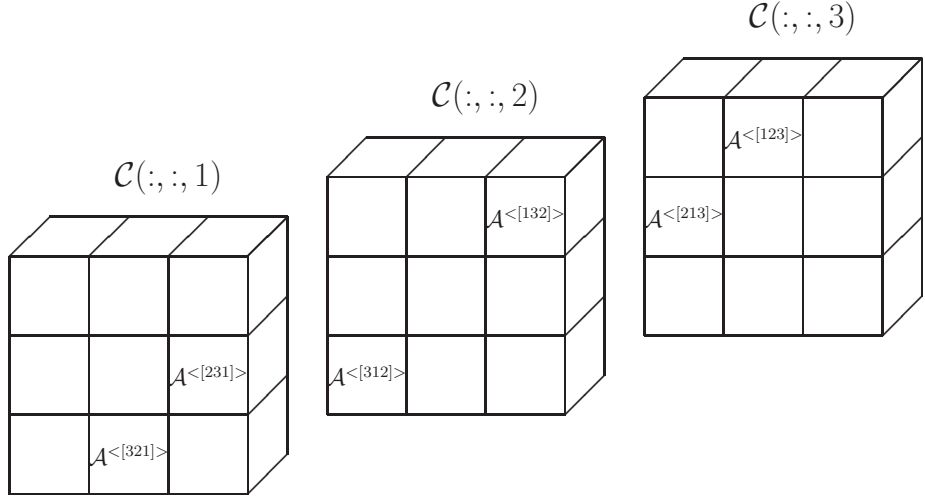

The last topic to cover in our order-3 preview is the generalization of the
Rayleigh quotient $\phi_{\cal A}^{(sym)}$ defined in (1.4). If
\inv{\cal A}{n_{1} \times n_{2} \times n_{3}}, ${\cal C} = \mbox{\bf sym}({\cal A})$,
$N = n_{1}+n_{2}+n_{3}$, and \inv{x}{N}, then $\phi_{\cal A}^{(sym)}$ is defined by
\begin{equation}
\phi_{\cal A}^{(sym)}(x) \;=\; \frac{1}{3!}
\left(
\sum_{i_{1}=1}^{N} \sum_{i_{2}=1}^{N} \sum_{i_{3}=1}^{N}
{\cal C}(i_{1},i_{2},i_{3}) x(i_{1})x(i_{2})x(i_{3})
\right) \Bigg/
\norm{x}_{2}^{3}
\end{equation}
It will be shown in section 4.3 that if
\[
x \;=\; \left[ \begin{array}{c}u\\v\\w \end{array} \right]\!\!
 \begin{array}{l}
\mbox{\small$\} n_{1}$} \\
\mbox{\small$\} n_{2}$} \\
\mbox{\small$\} n_{3}$} 
\end{array}
\]
satisfies $\nabla \phi_{\cal A}^{(sym)}(x) = 0$, then 
\[
\nabla_{\!\!u}\,\phi_{\cal A}(u,v,w) \:=\: 0\qquad
\nabla_{\!\!v}\,\phi_{\cal A}(u,v,w) \:=\: 0\qquad
\nabla_{\!\!w}\,\phi_{\cal A}(u,v,w) \:=\: 0
\]
where $\nabla_{\!\!z}$ refers to the gradient
with respect to the components in vector $z$.
Moreover, it will be shown that  
\[
x_{\mbox{{\tiny{+$-$}}}} = \left[ \begin{array}{r} \!\!u\!\! \\ \!\!v \!\!\\ \!\! -\!w\!\! \end{array} \right] \qquad
x_{\mbox{\tiny{$-+$}}} = \left[ \begin{array}{r} u \\ -v \\ w \end{array} \right] \qquad
x_{\mbox{\tiny{$--$}}} = \left[ \begin{array}{r} u \\ -v \\ -w \end{array} \right]
\] 
are also stationary vectors 
for $\phi_{\cal A}^{(sym)}$  and
\[
\phi_{\cal A}(u,v,w) \;=\;
\phi_{\cal A}^{(sym)}(x) \;=\; \phi_{\cal A}^{(sym)}(x_{\mbox{\tiny{$--$}}}) \;=\;
-\phi_{\cal A}^{(sym)}(x_{\mbox{\tiny{$-+$}}}) \;=\; -\phi_{\cal A}^{(sym)}(x_{\mbox{\tiny{$+-$}}}).
\]

\section{The Symmetric Embedding}

Block matrix manipulation is such a fixture in numerical linear algebra that
we take for granted the correctness of facts like
\begin{equation}
\left[ \begin{array}{cc} A_{11} & A_{12} \\ A_{21} & A_{22}\rule{0pt}{15pt} \end{array} \right]^{T} \;=\;
\left[ \begin{array}{cc} A_{11}^{T} & A_{21}^{T} \\ A_{12}^{T} & A_{22}^{T} \rule{0pt}{15pt}\end{array} \right]. \label{blocktrans}
\end{equation}
Formal verification requires showing that the $(i,j)$ entries on both sides of the equation are
equal for all valid $ij$ pairs.

The symmetric embedding of a tensor involves generalizations of both transposition and blocking
so this section begins by discussing these notions and establishing the tensor analog of
(2.1). 
Since vectors of subscripts are prominent in the presentation, we elevate their
notational status with boldface font, e.g., $\mathbf{p} \;=\; [\: 4\:1 \:2 \:3\:]$. We let
$\mathbf{1}$ denote the vector of ones and assume that dimension is clear from context. More
generally, if
$N$ is an integer, then $\mathbf{N}$ is the vector of all $N$'s. Finally, if
$\mathbf{i}$ and $\mathbf{j}$ have equal length, then $\mathbf{i} \leq \mathbf{j}$ means
that $i_{k} \leq j_{k}$ for all $k$.

\subsection{Blocking}

If $s$ and  $t$ are integers with $s \leq t$, then (as in {\sc Matlab}) let
$s:t$ denote  the row vector $[s,s+1, \cdots ,t]$. We refer to a vector with
this form as an {\em index range vector}. The act of blocking an $m_{1}$-by-$m_{2}$ matrix $C$ is the act of partitioning 
the index range vectors $1:m_{1}$ and $1:m_{2}$: 
\begin{equation}
\mathbf{r}^{(1)} \;= \;1:m_{1} \;=\; \left[ \begin{array}{c|c|c} \mathbf{r}_{1}^{(1)} & \cdots & \mathbf{r}_{b_{1}}^{(1)} \end{array} \right]
\qquad
\mathbf{r}^{(2)} \;= \;
1:m_{2} \;=\; \left[ \begin{array}{c|c|c} \mathbf{r}_{1}^{(2)} & \cdots & \mathbf{r}_{b_{2}}^{(2)} \end{array} \right] 
\end{equation}
Given  (2.2), we
are able to regard $C$ as a $b_{1} \! \times\! b_{2}$ block matrix $\left( C_{i_{1},i_{2}} \right)$
where block $C_{i_{1},i_{2}}$ has 
$\mbox{\tt length}(\mathbf{r}_{i_{1}}^{(1)})$ rows and $\mbox{\tt length}(\mathbf{r}_{i_{2}}^{(2)})$ columns.
It is easy (although messy) to ``locate'' a particular entry of a particular block. Indeed,
\[
C_{i_{1},i_{2}}(j_{1},j_{2}) \;=\; C(\:\rho_{i_{1}}^{(1)}+j_{1}\:,\:\rho_{i_{2}}^{(2)}+j_{2}\:)
\]
where
\begin{equation}\label{rhodef}
\rho_{i_{k}}^{(k)} \;=\; \mbox{\tt length}(\mathbf{r}_{1}^{(k)}) + \mbox{\tt length}(\mathbf{r}_{2}^{(k)}) + \cdots + 
\mbox{\tt length}(\mathbf{r}_{i_{k}-1}^{(k)})
\end{equation}
 for $k=1:2.$

To block  an order-$d$ tensor $\inTens{C}{}{m_{1}}{m_{d}}$ we proceed analogously.  
The index-range vectors
$1:m_{1},\ldots,1:m_{d}$ are partitioned
\begin{equation}
\mathbf{r}^{(k)} \;=\; 1:m_{k} \;=\;\left[ \begin{array}{c|c|c} \mathbf{r}^{(k)}_{1} & \cdots & 
\mathbf{r}^{(k)}_{b_{k}} \end{array} \right] \qquad k=1:d\label{indexvpart}
\end{equation} 
and this permits us to regard $\cal C$ as a $b_{1} \! \times \! \cdots \!\times\! b_{d}$ 
block tensor. If $\mathbf{i} = [i_{1},\ldots,i_{d}]$, then
the $\mathbf{i}$-th block is the subtensor
\[
{\cal C}_{\mathbf{i}} \;=\; {\cal C}_{i_{1},\ldots,i_{d}} \;=\; {\cal C}(\mathbf{r}_{i_{1}}^{(1)},\ldots,\mathbf{r}_{i_{d}}^{(d)}).
\]
If  $\mathbf{j} = [j_{1},\ldots,j_{d}]$, then the $\mathbf{j}$-th entry of this subtensor
is given by
\begin{equation}
\mathcal{C}_{\mathbf{i}}(\mathbf{j}) \;=\; \mathcal{C}(\:\rho_{i_{1}}^{(1)}+j_{1}\:,\ldots,\:\rho_{i_{d}}^{(d)}+j_{d}) \in \mathbb{R}
\end{equation}
where $\rho_{i_{k}}^{(k)} $ is specified by (2.3) for $k=1:d$.

To illustrate equations (2.3)-(2.5), if
 \inT{C}{9}{7}{5}{6} and
\[
\begin{array}{lclcl}
1:9 &=& \left[ \begin{array}{cccc|cc|ccc} 1 \!&\! 2 \!&\! 3 \!&\! 4 & 5 \!&\! 6 & 7 \!&\! 8 \!&\! 9 \end{array} \right] & \;& (b_{1} = 3)\\
1:7 &=& \left[ \begin{array}{ccccc|cc} 1 \!&\! 2 \!&\! 3 \!&\! 4 \!&\! 5 & 6 \!&\! 7  \end{array} \right] & \;& (b_{2} = 2) \rule{0pt}{12pt}\\
1:5 &=& \left[ \begin{array}{cccc|c} 1 \!&\! 2 \!&\! 3 \!&\! 4 & 5  \end{array} \right]& \;& (b_{3} = 2) \rule{0pt}{12pt}\\
1:6 &=& \left[ \begin{array}{cc|cc|cc} 1 \!&\! 2 & 3 \!&\! 4 & 5 \!&\! 6 \end{array} \right]& \;& (b_{4} = 3) \rule{0pt}{12pt}
\end{array}\; ,
\]
then we are choosing to regard
 $\cal C$ as a $3\times 2 \times 2 \times 3$ block tensor. Thus, if
$\mathbf{i} = [3\:1\:2\:1]$ then
${\cal C}_{\mathbf{i}} = {\cal C}(7:9,1:5,5:5,1:2)$ and
\[
{\cal C}_{\mathbf{i}}(\mathbf{j}) \;=\;
{\cal C}(6+j_{1}\:,\:j_{2}\:,\:4+j_{3}\:,\:j_{4})
\]
where $\mathbf{1} \leq \mathbf{j} \leq [3\:5\:1\:2]$.

\subsection{Tensor Transposition}

If 
$ \inTens{A}{}{n_{1}}{n_{d}}$ and $\mathbf{p} = [p_{1},\ldots,p_{d}]$
is a permutation of $1:d$, then $\inTens{A^{<\mathbf{p}>}}{}{n_{p_{1}}}{n_{p_{d}}}$ denotes the
$\mathbf{p}$-transpose of $\cal A$ defined by
\[
{\cal A}^{<\mathbf{p}>}(j_{p_{1}},\ldots,j_{p_{d}}) \;=\; {\cal A}(j_{1},\ldots,j_{d}) 
\]
where $1\leq j_{k} \leq n_{k}$ for $k=1:d.$ A more succinct way of saying the
same thing is
\[
{\cal A}^{<\mathbf{p}>}(\mathbf{j}(\mathbf{p})) = {\cal A}(\mathbf{j}) \qquad \mathbf{1} \leq \mathbf{j} \leq \mathbf{n}.
\]
 If $\cal A$ is an order-2 tensor,
then ${\cal A}^{<\,[2\;1]\,>}(j_{2},j_{1}) = {\cal A}(j_{1},j_{2})$.  It is also easy to verify that if $\mathbf{f}$ and $\mathbf{g}$ are both permutations of $1:d$, then
\begin{equation}
( {\cal A}^{<\mathbf{f}>})^{<\mathbf{g}>} \;=\; {\cal A}^{<\mathbf{f}(\mathbf{g})>}.
\end{equation}

A transposition of a block tensor renders another block tensor. The following lemma makes this
precise and generalizes (2.1).

\medskip

\begin{lemma} Suppose $\inTens{C}{}{m_{1}}{m_{d}}$ is a $b_{1}\times \cdots \times b_{d}$ block tensor with block dimensions defined by the
partitioning (2.4). Let ${\cal C}_{\mathbf{i}}$ denote its $\mathbf{i}$-th block
 where $\mathbf{i} = [i_{1},\ldots,i_{d}]$.
If $\mathbf{p} = [\:p_{1},\ldots,p_{d}\:]$ is a permutation of $1:d$ and
${\cal B} = {\cal C}^{<\mathbf{p}>}$, then the tensor $\inTens{B}{}{m_{p_{1}}}{m_{p_{d}}}$ is a 
$b_{p_{1}} \times \cdots \times b_{p_{d}}$ block tensor where each block 
${\cal B}_{\mathbf{i}(\mathbf{p})}$ is
defined by ${\cal B}_{\mathbf{i}(\mathbf{p})} = {\cal C}_{\mathbf{i}}^{<\mathbf{p}>}.$
\end{lemma}

\medskip

\begin{proof}
If $1 \leq j_{k} \leq m_{k}$ for $k=1:d$, then from (2.4) and (2.5) we have
\[
{\cal C}_{\bf i}^{<{\bf p}>}(j_{p_{1}},\ldots,j_{p_{d}}) \;=\; {\cal C}_{\bf i}(j_{1},\ldots,j_{p}) 
\;=\; {\cal C}(\rho_{i_{1}}^{(1)} + j_{1},\ldots, \rho_{i_{d}}^{(d)} + j_{d}) 
\]
On the other hand, ${\cal B} = {\cal C}^{<{\bf p}>}$ and so
\[
{\cal C}(\rho_{i_{1}}^{(1)} + j_{1},\ldots, \rho_{i_{d}}^{(d)} + j_{d}) 
\;=\; {\cal B}(\rho_{i_{p_{1}}}^{(p_{1})} + j_{p_{1}},\ldots, \rho_{i_{p_{d}}}^{(p_{d})} + j_{p_{d}}) \;=\; {\cal B}_{\mathbf{i}(\mathbf{p})}(j_{p_{1}},\ldots,j_{p_{d}}).
\]
Thus, ${\cal B}_{\mathbf{i}(\mathbf{p})}(\mathbf{j}(\mathbf{p}))\;=\;
{\cal C}_{\mathbf i}^{<{\bf p}>}(\mathbf{j}(\mathbf{p}))$ for all $\mathbf{j}$ , i.e.,
${\cal B}_{\mathbf{i}({\bf p})} = {\cal C}_{\mathbf{i}}^{<{\bf p}>}.$
\end{proof}

\subsection{The sym$(\cdot)$ Operation}
An  order-$d$ tensor \inTens{C}{}{N}{N} is {\em symmetric} if 
${\cal C} = {\cal C}^{<\mathbf{p}>}$
for any permutation $\mathbf{p}$ of $1:d$. The tensor
analog of (1.1) involves constructing an order-$d$ symmetric tensor 
$\mathbf{sym}(\cal A)$ whose blocks are
either zero or carefully chosen transposes of $\cal A$. In particular,
if $\inTens{A}{}{n_{1}}{n_{d}}$, then
\[
\inTens{\mathbf{sym}({A})}{}{N}{N} \qquad N = n_{1}+\cdots n_{d}
\]
is a block tensor  defined by the partitioning
$1:N = [ \;\mathbf{r}_{1} \;|\; \cdots \;|\; \mathbf{r}_{d}\;]$ 
where 
\begin{equation}
\mathbf{r}_{k} = (1+n_{1}+\cdots + n_{k-1}):(n_{1}+\cdots + n_{k})\qquad k=1:d.
\end{equation}
The $\mathbf{i}$-th block of ${\cal C} = \mathbf{sym}({\cal A})$ is given by
\[
{\cal C}_{\mathbf{i}} \;=\;
\left\{ 
\begin{array}{ll}
A^{<\,\mathbf{i}\,>} & \mbox{if $\mathbf{i}$ is a permutation of $1:d$}\\
\\
0 & \mbox{otherwise} 
\end{array}
\right.
\]
for all $\mathbf{i}$ that satisfy $\mathbf{1} \leq \mathbf{i} \leq \mathbf{d}$.
Note that ${\cal C}_{\mathbf{i}}$ is $n_{i_{1}}\!\times \! n_{i_{2}} \! \times \! \cdots \!\times \!n_{i_{d}}$.
We confirm that
$\mathbf{sym}({\cal A})$ is symmetric.

\medskip

\begin{lemma} If $\inTens{A}{}{n_{1}}{n_{d}}$ and ${\cal C} = \mathbf{sym}({\cal A})$, then
$\cal C$ is symmetric.
\end{lemma}

\medskip

\begin{proof}
Let $\mathbf{p}$ be an arbitrary permutation of $1:d$. We must show that if ${\cal B} = 
C^{<\,\mathbf{p}\,>}$ then ${\cal B} = {\cal C}$. Since $\cal C$ as a block tensor is
$d \! \times \! d \! \times \!\cdots \!\times \! d$, it follows from Lemma 2.1 that
${\cal B}$ has the same block structure and
\[
{\cal B}_{\mathbf{i}(\mathbf{p})} \;=\; {\cal C}_{\mathbf{i}}^{<\,\mathbf{p}\,>}
\]
for all $\mathbf{i}$ that satisfy $\mathbf{1} \leq \mathbf{i} \leq \mathbf{d}$. 
If $\mathbf{i}$ is a permutation of $1:d$, then
${\cal C}_{\mathbf{i}} = A^{<\,\mathbf{i}\,>}$ and by using (2.6) we conclude that
\[
{\cal B}_{\mathbf{i}(\mathbf{p})} \;=\;  ( {\cal A}^{<\mathbf{i}>})^{<\mathbf{p}>} \;=\; {\cal A}^{<\mathbf{i}(\mathbf{p})>} \;=\; {\cal C}_{\mathbf{i}(\mathbf{p})}
\]
If $\mathbf{i}$ is {\em not} a permutation of $1:d$, then both ${\cal C}_{\mathbf{i}}$ and ${\cal C}_{\mathbf{i}(\mathbf{p})}$
are zero and so
\[
{\cal B}_{\mathbf{i}(\mathbf{p})} = {\cal C}_{\mathbf{i}}^{<\,\mathbf{p}\,>} = 0 = {\cal C}_{\mathbf{i}(\mathbf{p})}.
\]
Since ${\cal B}$ and ${\cal C}$ agree block-by-block, they are the same.
\end{proof}

\section{Orderings, Unfoldings, and Summations}

In numerical multilinear algebra it is frequently necessary to reshape a given tensor into a vector 
or a matrix and vice versa. In this section we collect results that make these maneuvers precise.

\subsection{The col Ordering}
If $\mathbf{i}$ and $\mathbf{s}$ are
length-$e$ index vectors and $\mathbf{1} \leq \mathbf{i} \leq \mathbf{s}$, then we define the integer-valued function
{\em ivec}  by
\begin{eqnarray*}
\mbox{\em ivec}(\mathbf{i},\mathbf{s}) &\;=\;& 
i_{1} + 
(i_{2}-1)s_{1} + 
\cdots + 
(i_{e}-1)(s_{1}\cdots s_{e-1})
\end{eqnarray*}
If $\inTens{F}{}{s_{1}}{s_{e}}$, then 
\inv{v = \mbox{\bf vec}({\cal F})}{s_{1}\cdots s_{e}} is the column vector
defined by
\[
v(\mbox{\em ivec}(\mathbf{i},\mathbf{s})) = {\cal F}(\mathbf{i}) \qquad 
\mathbf{1} \leq \mathbf{i} \leq \mathbf{s}.
\]
Note that if $e = 2$, then $\cal F$ is a matrix and {\bf vec}$({\cal F})$ stacks its columns.
We also observe that if $\inv{w_{k}}{s_{k}}$ for $k=1:e$, then
\begin{equation}
w = w_{e}\T \cdots \T w_{1} \qquad \Leftrightarrow \qquad w(\mbox{\em ivec}(\mathbf{i},\mathbf{s})) = w_{1}(i_{1})\cdots w_{e}(i_{e}).
\end{equation}

\subsection{Modal Unfoldings}
In the gradient calculations that follow, it is particularly convenient to
``flatten'' the given tensor $\inTens{A}{}{n_{1}}{n_{d}}$ into a matrix.
If 
\begin{eqnarray}
\tilde{\mathbf{n}} &=& [\:\mathbf{n}(1:k-1) \;\;\; \mathbf{n}(k+1:d)\:],\\
\tilde{\mathbf{i}} &=& [\:\:\mathbf{i}(1:k-1) \; \;\;\; \mathbf{i}(k+1:d)\:], \rule{0pt}{13pt}
\end{eqnarray}
 then the {\em mode-$k$ unfolding}
${\cal A}_{(k)}$ is defined 
by
\begin{equation} \label{modekunfolding}
{\cal A}_{(k)}(i_{k}, \mbox{\em ivec}(\tilde{\mathbf{i}},\tilde{\mathbf{n}}) ) 
\;=\; {\cal A}(\mathbf{i}) \qquad \mathbf{1} \leq \mathbf{i} \leq \mathbf{n}.
\end{equation}
This matrix has $n_k$ rows and $n_1\cdots n_{k-1}n_{k+1}\cdots n_d$ columns. A third-order instance of this important concept is displayed in equation (1.6). 
We mention that there are other ways to order the columns in ${\cal A}_{(k)}$.
See \cite{koldaToolbox}. 

While the columns of ${\cal A}_{(k)}$ are mode-$k$ fibers, its rows are reshapings
of its {\em  mode-$k$  subtensors}. In particular, if $1\leq r \leq n_{k}$, then
\[
{\cal A}_{(k)}(r,:) \;=\; \mbox{\bf vec}({\cal B}^{(r)})^{T}
\]
where the mode-$k$ subtensor ${\cal B}^{(r)}$ has order $d-1$ and is defined by
\[
{\cal B}^{(r)}(i_{1},\ldots,i_{k-1},i_{k+1},\ldots,i_{d}) \;=\;
{\cal A}(i_{1},\ldots,i_{k-1},r,i_{k+1},\ldots,i_{d}).
\]
The partitioning of an order-$d$ tensor into order-$(d-1)$ tensors is just a generalization of
partitioning a matrix into its columns.
\medskip

\subsection{Summations}

It is handy to have a multi-index summation notation in order to describe general versions of the summations that appear in (1.5) and (1.12). If $\mathbf{n}$ is a length-$d$
index vector, then
\[
\sum_{\mathbf{i}=\mathbf{1}}^{\mathbf{n}} \;\equiv\; 
\sum_{i_{1}=1}^{n_{1}} \cdots \sum_{i_{d}=1}^{n_{d}}.
\]
The summation that defines the multilinear Rayleigh quotient (1.5) can
be written in matrix-vector terms.

\medskip

\begin{lemma}
If $\inTens{A}{}{n_{1}}{n_{d}}$ and \inv{u_{k}}{n_{k}} for $k=1:d$, then
\begin{equation}
\sum_{\mathbf{i} = \mathbf{1}}^{\mathbf{n}} {\cal A}(\mathbf{i})u_{1}(i_{1})\cdots u_{d}(i_{d}) \;=\; \mathbf{vec}({\cal A})^{T} u_{d} \T \cdots \T u_{1}.
\end{equation}
Moreover, for $k=1:d$ we have
\begin{equation}
\sum_{\mathbf{i}=1}^{\mathbf{n}} {\cal A}(\mathbf{i})u_{1}(i_{1}) \cdots u_{d}(i_{d})
\;=\; u_{k}^{T}\, {\cal A}_{(k)} \,\tilde{u}_{k}
\end{equation}
where
\begin{equation}
\tilde{u}_{k} \;=\; (u_{d} \T \cdots \T u_{k+1} \T u_{k-1} \T \cdots \T u_{1}).
\end{equation}
\end{lemma}
\begin{proof}
If $a =  \mathbf{vec}({\cal A})$ and $b = u_{d} \T \cdots \T u_{1}$, then
using the definition of {\bf vec} and equations (3.1)-(3.4),
we have
\[
\sum_{\mathbf{i} = \mathbf{1}}^{\mathbf{n}} {\cal A}(\mathbf{i})u_{1}(i_{1})\cdots u_{d}(i_{d})
\;=\;
\sum_{\mathbf{i} = \mathbf{1}}^{\mathbf{n}} 
a(\mbox{\em ivec}(\mathbf{i},\mathbf{n}))\cdot b(\mbox{\em ivec}(\mathbf{i},\mathbf{n}))
\;=\; \sum_{\mathbf{i} = \mathbf{1}}^{\mathbf{n}} a(\mathbf{i})b(\mathbf{i}) \;=\; a^{T}b.
\]
This proves (3.5).
Using the modal subtensor interpretation of
${\cal A}_{(k)}$ that we discussed in \S3.2 and definitions (3.2) and (3.3), we have
\begin{eqnarray*}
\sum_{\mathbf{i}=1}^{\mathbf{n}} {\cal A}(\mathbf{i})u_{1}(i_{1}) \cdots u_{d}(i_{d})
&=& \sum_{i_{k}=1}^{n_{k}}u_{k}(i_{k})\left( \sum_{\tilde{\mathbf{i}}=\mathbf{1}}^{\tilde{\mathbf{n}}}
{\cal B}^{(i_{k})}(\tilde{\mathbf{i}})\tilde{u}(\tilde{\mathbf{i}})\right)\\
&=&
 \sum_{i_{k}=1}^{n_{k}}u_{k}(i_{k})\left( {\cal A}_{(k)}(i_{k},:)\tilde{u}_{k} \right)
\;=\; u_{k}^{T} {\cal A}_{(k)} \tilde{u}_{k}
\end{eqnarray*}
\vspace*{-.1in}
which establishes (3.6).
\end{proof}

\bigskip
Summations that involve symmetric tensors are important in later sections. The following notation
for the multiple Kronecker product of a single vector is handy:

\vspace*{-.1in}
\[
x^{\otimes d} = \underbrace{x\;\otimes\; \cdots\; \otimes\; x}_{d \mbox{ times}}.
\]
Note that if \inv{x}{N}, then \inv{x^{\otimes d}}{N^{d}}.

\medskip

\begin{lemma}
If $\inTens{C}{}{N}{N}$ is a symmetric order-$d$ tensor and \inv{x}{N}, then
\begin{equation}
\sum_{\mathbf{i}=1}^{\mathbf{N}} {\cal C}(\mathbf{i})x(i_{1}) \cdots x(i_{d})
\;=\; x^{T}\, {\cal C}_{(1)} \, x^{\otimes (d-1)}
\end{equation}
\end{lemma}
\vspace*{-.1in}
\begin{proof}
This follows from Lemma 3.1 by setting $n_{k} = N$ and $u_{k} = x$ for $k=1:d$. Note that
because ${\cal C}$ is symmetric, ${\cal C}_{(1)} = \cdots = {\cal C}_{(d)}$.
\end{proof}

\medskip

The summation (3.8) has a special characterization if ${\cal C} = \mathbf{sym}({\cal A})$. To pursue this we will have to navigate $\cal C$'s block structure and to that end we define the index
vectors $\mathbf{L}$ and $\mathbf{R}$ as follows:
\begin{equation}\label{LRvecs}
\mathbf{L} \;=\; \left[ \begin{array}{c} 1 \\ n_{1}+1 \\ \vdots \\ n_{1}+\cdots+n_{d-1} +1
\end{array} \right]
\qquad
\mathbf{R} \;=\; \left[ \begin{array}{c} n_{1} \\ n_1 + n_{2} \\ \vdots \\ n_{1}+\cdots+n_{d}
\end{array} \right].
\end{equation}
Note that if $\mathbf{1} \leq \mathbf{p} \leq \mathbf{d}$, then
\[
{\cal C}_{\mathbf p} \;=\; {\cal C}(\mathbf{L}(p_{1}):\mathbf{R}(p_{1}),\ldots,
\mathbf{L}(p_{d}):\mathbf{R}(p_{d}))
\]
 is ${\cal C}$'s $\mathbf{p}$-th block.

\medskip

\begin{lemma}
Suppose $\inTens{A}{}{n_{1}}{n_{d}}$, ${\cal C} \;=\; \mbox{\bf sym}({\cal A})$, and $N = n_{1}+\cdots + n_{d}$. If \inv{x}{N} is partitioned as follows
\vspace*{-.1in}
\[
x \;=\; \left[ \begin{array}{c} u_{1} \\ \vdots \\ u_{d} \end{array} \right]
\qquad \inv{u_{k}}{n_{k}},
\]
and $\tilde{u}_{1},\ldots,\tilde{u}_{d}$ are defined by (3.7), then

\vspace*{-.2in}
\begin{equation}
{\cal C}_{(1)} x^{\otimes (d-1)} \;=\; (d-1)! \left[ \begin{array}{c}
{\cal A}_{(1)} \tilde{u}_{1} \\
\vdots\\
{\cal A}_{(d)} \tilde{u}_{d} \end{array} \right]
\end{equation}
and
\begin{equation}
\sum_{\mathbf{j}=\mathbf{1}}^{\mathbf{N}} {\cal C}(\mathbf{j})x(j_{1})\cdots x(j_{d})
\;=\; d! \sum_{\mathbf{i}=\mathbf{1}}^{\mathbf{n}} {\cal A}(\mathbf{i}) u_{1}(i_{1})\cdots u_{d}(i_{d}).
\end{equation}
\end{lemma}
\begin{proof}
If $v = {\cal C}_{(1)}x^{\otimes (d-1)}$ and 
\begin{equation}
e_j = I_{N}(:,j) \;=\;
\left[ \begin{array}{c} w_{1} \\ \vdots \\ w_{d} \end{array} \right]
\end{equation}
is partitioned conformally with $x$, then for $j=1:N$ we have
\begin{eqnarray*}
v(j) &=&\sum_{\mathbf{i}(2:d)=\mathbf{1}}^{\mathbf{N}} \mathcal{C}(j,i_2,\ldots,i_d) x(i_2)\cdots x(i_d) \\
&=&
 \sum_{\mathbf{i}=\mathbf{1}}^\mathbf{N}\mathcal{C}(\mathbf{i}) e_j(i_1)x(i_2) \cdots x(i_d)\\
&=& \sum_{\mathbf{p}= \mathbf{1}}^{\mathbf{d}}
\sum_{\mathbf{i}=\mathbf{L}(\mathbf{p})}^{\mathbf{R}(\mathbf{p})}
\mathcal{C}(\mathbf{i}) e_j(i_1)x(i_2) \cdots x(i_d)\\
&=& \sum_{\mathbf{p}= \mathbf{1}}^{\mathbf{d}}\left(
\sum_{\mathbf{k}=\mathbf{1}}^{\mathbf{n}(\mathbf{p})}
\mathcal{C}_{\mathbf{p}}(\mathbf{k}) 
w_{p_{1}}(k_{1})u_{p_{2}}(k_2) \cdots u_{p_{d}}(k_d)\right)
\end{eqnarray*}
Now suppose that $\mathbf{L}(q) \leq j \leq \mathbf{R}(q)$, $j = \mathbf{L}(q) + r -1$. From (3.10) we 
must show that $v_{j}$ is the $r$th component of ${\cal A}_{(q)} \tilde{u}_{q}$.

To that end observe that $\mathcal{C}_{\mathbf{p}}(\mathbf{k}) w_{p_{1}}(k_{1})$  is necessarily zero
unless 
$p_{1} = q$, $k_{1}=  r$, and $\mathbf{p}$ is  a permutation of $1:d$. 
Assuming this to be the case and defining the vectors $v_{1},\ldots,v_{d}$ by
\[
v_{i} \;=\; \left\{ \begin{array}{ll} u_{i} & \mbox{if $i\neq q$} \\ w_{q} & \mbox{otherwise}
\rule{0pt}{16pt}
\end{array}\right.,
\]
we see using (3.6) that

\begin{eqnarray*}
\sum_{\mathbf{k}=\mathbf{1}}^{\mathbf{n}(\mathbf{p})}
\mathcal{C}_{\mathbf{p}}(\mathbf{k}) 
w_{p_{1}}(k_{1})u_{p_{2}}(k_2) \cdots u_{p_{d}}(k_d) &=&
\sum_{\mathbf{k}=\mathbf{1}}^{\mathbf{n}(\mathbf{p})}
{\cal A}^{<p>}(\mathbf{k})v_{p_{1}}(k_{1})v_{p_{2}}(k_2) \cdots v_{p_{d}}(k_d)\\
&=&
\sum_{\mathbf{k}=\mathbf{1}}^{\mathbf{n}}
{\cal A}(\mathbf{k})v_{1}(k_{1})v_{2}(k_2) \cdots v_{d}(k_d)\\
&=& v_{q}^{T} {\cal A}_{(q)} v_{d} \T \cdots \T v_{q+1} \T v_{q-1} \T \cdots \T v_{1} \rule{0pt}{18pt} \\
&=& w_{q}^{T}  {\cal A}_{(q)} u_{d} \T \cdots \T u_{q+1} \T u_{q-1} \T \cdots \T u_{1} \rule{0pt}{18pt}\\
&=& w_{q}^{T} {\cal A}_{(q)} \tilde{u}_{q}. \rule{0pt}{18pt}
\end{eqnarray*}
Observe that the number of $\mathbf{p}$ that satisfy $\mathbf{1} \leq \mathbf{p} \leq \mathbf{d}$ 
subject to the constraint $p_{1} = q$ is $(d-1)!$ and conclude from (3.12) that $w_{q} = I_{n_{q}}(:,r)$. It follows 
that 
\[ 
v(j) \;=\; 
\sum_{\mathbf{p}\;=\; \mathbf{1}}^{\mathbf{d}} w_{q}^{T} {\cal A}_{(q)} \tilde{u}_{q}\;=\;
(d-1)!\left[ {\cal A}_{(q)} \tilde{u}_{q}\right]_{r}.
\]
This establishes (3.10). Equation (3.11) follows from
\[
x^{T} {\cal C}_{(1)}x^{\otimes (d-1)} \;=\;
 \sum_{k=1}^{d} (d-1)! u_{k}^{T}{\cal A}_{(k)}\tilde{u}_{k} 
\]
and Lemmas 3.1 and 3.2.
\end{proof}

\section{Rayleigh Quotients and Stationary Values}

Suppose $\inTens{A}{}{n_{1}}{n_{d}}$ and \inv{u_{k}}{n_{k}} for $k=1:d$. Analogous to (1.3) we define the multilinear Rayleigh Quotient
\begin{equation}
\phi_{\cal A}(u_{1},\ldots,u_{d}) \;=\; 
\left( \sum_{\mathbf{i} = \mathbf{1}}^{\mathbf{n}} {\cal A}(\mathbf{i})u_{1}(i_{1})\cdots u_{d}(i_{d}) \right) \Bigg/\left( \norm{u_{1}}_{2} \cdots \norm{u_{d}}_{2} \right)
\end{equation}
If ${\cal C} = \mathbf{sym}({\cal A})$, $N = n_{1}+\cdots n_{d}$, and \inv{x}{N}, then corresponding to
(1.4) we have
\begin{equation}
\phi_{\cal A}^{(sym)}(x) \;=\; \frac{1}{d!}
\left( \sum_{\mathbf{i} = \mathbf{1}}^{\mathbf{N}} {\cal C}(\mathbf{i})x(i_{1})\cdots x(i_{d}) \right) \Bigg/ \left( \norm{x}_{2}\right)^{d} 
\end{equation}
In this section we examine these multilinear Rayleigh quotients, specify their gradients, and relate the singular values of $\cal A$ to the eigenvalues of $\mathbf{sym}({\cal A})$.

\subsection{The Singular Values of a General Tensor}

The gradient of $\phi_{\cal A}(u_{1},\ldots,u_{d})$ relates to a collection of matrix-vector
products that involve the modal unfoldings of $\cal A$ and Kronecker products of
the $u$-vectors.

\medskip

\begin{theorem}\label{gradphiA}
If $\inTens{A}{}{n_{1}}{n_{d}}$ and for $k=1:d$ the vectors \inv{u_{k}}{n_{k}} each have unit 2-norm, then
\[
\nabla \phi_{\cal A}(u_{1},\ldots,u_{d}) \;= \; 
\left[ 
\begin{array}{c}
{\cal A}_{(1)} \tilde{u}_{1}\\
\vdots\\
{\cal A}_{(d)} \tilde{u}_{d}
\end{array}\right]
\;-\; \phi_{\cal A}(u_{1},\ldots,u_{d})\cdot
\left[ \begin{array}{c} u_{1} \\ \vdots \\ u_{d}
\end{array} \right]
\]
where
$\tilde{u}_{k} \;=\; (u_{d} \T \cdots \T u_{k+1} \T u_{k-1} \T \cdots \T u_{1})$.
\end{theorem}
\begin{proof}
From Lemma 3.1 we have
\[
\phi_{\cal A}(u_{1},\ldots,u_{d}) 
\;=\; \left(u_{k}^{T}\, {\cal A}_{(k)} \,\tilde{u}_{k} \right)
/ \left(\norm{u_{1}}_{2} \cdots \norm{u_{d}}_{2} \right).
\]
For $k=1:d$ we have
$
\nabla_{u_{k}} \!\left(u_{k}^{T}\, {\cal A}_{(k)} \,\tilde{u}_{k} \right)\:=\:
{\cal A}_{(k)} \,\tilde{u}_{k}
$
and
$
\nabla_{u_{k}}\!\left( \norm{u_{1}}_{2} \cdots \norm{u_{d}}_{2}\right) \:=\: u_{k}.
$
and so
\begin{eqnarray*}
\nabla_{u_{k}} \phi_{\cal A} &=&
\frac{
(\norm{u_{1}}_{2}\cdots \norm{u_{d}}_{2})
{\cal A}_{(k)} \tilde{u}_{k} \:-\: (u_{k}^{T} {\cal A}_{(k)} \tilde{u}_{k}) u_{k}
}
{
(\norm{u_{1}}_{2}\cdots \norm{u_{d}}_{2})^{2}
}\\
&=&
{\cal A}_{(k)} \tilde{u}_{k}  \;-\; \phi_{\cal A}(u_{1},\ldots,u_{d}) u_{k}.\rule{0pt}{16pt}
\end{eqnarray*}
The theorem follows by simply ``stacking'' these subvectors of the  gradient.
\end{proof}

\medskip
The variational approach to tensor singular values and vectors set forth in \cite{lim2} 
is based on equating the gradient of $\phi_{\cal A}$ to zero.

\medskip

\begin{definition}The scalar $\sigma \in\RR$ is a \emph{singular value} of a general tensor 
 $\inTens{A}{}{n_1}{n_d}$ if there are unit vectors \inv{u_{k}}{n_{k}} such that
\begin{equation} \label{svaldef}
{\cal A}_{(k)} \tilde{u}_k = \sigma u_k,
\end{equation}
for  $k=1:d$. The vector $u_{k}$ is the {\em mode-$k$ singular vector} associated with $\sigma$.
\end{definition}

\medskip

\noindent
The normalization condition $u_k^T u_k =1$ is necessary, since if $v_k = a u_k$ for $k=1:d$ then ${\cal A}_{(k)} \tilde{v}_k = a^{d-1} {\cal A}_{(k)} \tilde{u}_k = a^{k-1} \sigma u_k = (a^{k-2}\sigma) v_k$ for any $a \in \RR$.
It can be shown that at least one singular value and associated singular vectors exist for any tensor (cf.~\cite{lim2}).

\subsection{The Eigenvalues of a Symmetric Tensor}

For a symmetric tensor ${\cal C}$, the stationary values  of $\phi_{\cal C}(x,\ldots,x)$ define the
notion of a tensor eigenvalue.
\medskip
\begin{theorem}
If $\inTens{\mathcal{C}}{}{N}{N}$ is symmetric and \inv{x}{N} has unit norm, then
\[
\nabla_x \phi_{\cal C}(x,\ldots,x)
= d \left({\cal C}_{(1)}x^{\otimes (d-1)} \:-\:\phi_{\cal C}(x,\ldots,x)  x \right).
\]
\end{theorem}
\begin{proof}
From Lemma 3.2 we have
\[
 \phi_{\cal C}(x,\ldots,x) \;=\; x^{T}{\cal  C}_{(1)} x^{\otimes(d-1)} /\, \norm{x}^{d}.
\]
Since
\[
\nabla_{x} x^{T}{\cal  C}_{(1)} x^{\otimes(d-1)} \;=\; d {\cal  C}_{(1)} x^{\otimes(d-1)}
\]
and
\[
\nabla_{x} (x^{T} x)^{d/2} \;=\; d(x^{T}x)^{d/2 - 1} x
\]
it follows that
\begin{eqnarray*}
\nabla_x \phi_{\cal C}(x,\ldots,x) &=& 
d \:\frac{(x^T x)^{d/2}  {\cal C}_{(1)} x^{\otimes (d-1)} \:- \:\left( x^T {\cal C}_{(1)} x^{ \otimes (d-1) }\right) (x^{T}x)^{d/2-1} x}{ (x^Tx)^d } \\ 
&=& d \left({\cal C}_{(1)}x^{\otimes (d-1)}\: - \:\left( x^T {\cal C}_{(k)} x^{\otimes (d-1)}\right)  x \right)\\
&=& d \left( {\cal C}_{(1)}x^{\otimes (d-1)} \: - \: \phi_{\cal C}(x,\ldots,x)  x\right)
\end{eqnarray*}
completing the proof of the theorem.
\end{proof}

\medskip
\noindent
By setting the gradient of $\phi_{\cal C}(x,\ldots,x)$ to zero we arrive at the notion
of a tensor eigenvalue \cite{qi4}.

\medskip

\begin{definition} 
If $\inTens{\cal C}{}{N}{N}$ is symmetric and \inv{x}{N} is a unit vector such that
\begin{equation}\label{evaldef}
{\cal C}_{(1)}x^{\otimes (d-1)} \: = \:\lambda x
\end{equation}
then $\lambda = \phi_{\cal C}(x,\ldots,x)$ is an eigenvalue of $\cal C$ and $x$ the
associated eigenvector.
\end{definition}

\medskip

\noindent
Note that if ${\cal C}_{(1)}x^{\otimes (d-1)} \: = \:\lambda x$ and \ins{\alpha},
then
${\cal C}_{(1)}(\alpha x)^{\otimes (d-1)} \: = \:(\alpha^{d-2} \lambda) (\alpha x)$.
Thus, we resolve a uniqueness issue by requiring tensor eigenvectors to have unit length,
something that is not necessary in the matrix ($d=2$) case.

In \cite{qi4,qi5} it is shown that  eigenvalues and associated eigenvectors always exist for symmetric tensors. Recently it has  been shown that a symmetric tensor has at most $((d-1)^N -1)/(d-2)$ eigenvalues, counted with multiplicity \cite{sturmfels}.

\subsection{The Eigenvalues of $\mathbf{sym}({\cal A})$}

Since ${\cal C} = \mathbf{sym}({\cal A})$ is so structured, we anticipate that the eigenvalue-defining equation $\nabla \phi_{\cal A}^{(sym)}(x) = 0$ will have some special features.
From the definitions (4.1) and (4.2) and Theorem 4.3, we have
\begin{equation}
\nabla \phi_{\cal A}^{(sym)}(x) \;=\; \frac{1}{d!} \nabla \phi_{\cal C}(x,\ldots,x)
\;=\; \frac{1}{(d-1)!}
 \left(
{\cal C}_{(1)} x^{\otimes(d-1)} \;-\; \phi_{\cal C}(x,\ldots,x)x
\right)
\end{equation}
We first characterize the gradient of $\phi_{\cal A}^{(sym)}$
in terms of matrix-vector products that involve ${\cal A}$'s modal unfoldings.

\medskip

\begin{theorem}
If $\inTens{\cal A}{}{n_{1}}{n_{d}}$ and $x$ has unit 2-norm, then
\begin{equation}\label{phisymgrad}
\nabla_{x} \phi_{\cal A}^{(sym)}(x) \;=\; 
\left[
\begin{array}{c}
{\cal A}_{(1)} \tilde{u}_{1} \\ \vdots \\ {\cal A}_{(d)} \tilde{u}_{d}
\end{array} \right]
\;-\;
d
\left[ \begin{array}{c}
(u_{1}^{T}{\cal A}_{(1)} \tilde{u}_{1}) u_{1} \\ \vdots \\
(u_{d}^{T}{\cal A}_{(d)} \tilde{u}_{d}) u_{d}
\end{array}
\right].
\end{equation}
\end{theorem}

\begin{proof}
From Lemmas 3.1 and 3.3 and the definitions (4.1) and 4.2) we have
\begin{eqnarray*}
\phi_{\cal A}^{(sym)}(x) &=& \frac{1}{d!}
\left( \sum_{\mathbf{i} = \mathbf{1}}^{\mathbf{N}} {\cal C}(\mathbf{i})x(i_{1})\cdots x(i_{d}) \right) \Bigg/ \left( \norm{x}_{2}\right)^{d} \\
&=& \left( \sum_{\mathbf{i} = \mathbf{1}}^{\mathbf{n}} {\cal A}(\mathbf{i})u_{1}(i_{1})\cdots u_(i_{d}) \right) \Bigg/ \left( \norm{x}_{2}\right)^{d} \\
&=& \frac{u_{k}^{T} {\cal A}_{(k)} \tilde{u}_{k}}{(u_{1}^{T}u_{1} + \cdots + u_{d}^{T}u_{d})^{d/2}}
\end{eqnarray*}
for $k=1:d$.
Since
\[
\nabla_{u_{k}}(u_{1}^{T}u_{1} \:+\: \cdots +\: u_{d}^{T}u_{d})^{d/2} \;=\;
d\cdot (u_{1}^{T}u_{1} \:+\: \cdots +\: u_{d}^{T}u_{d})^{d/2 - 1}u_{k} \;=\; d \cdot  u_{k}.
\]
Since $x^{T}x = u_{1}^{T}u_{1} + \cdots u_{d}^{T}u_{d}$,
we can conclude  that
\begin{eqnarray*}
\nabla_{u_{k}} \phi_{\cal A}^{(sym)}(x)  &=&
\frac{(x^{T}x)^{d/2} {\cal A}_{(k)} \tilde{u}_{k} - d\cdot (u_{k}^{T}{\cal A}_{(k)} \tilde{u}_{k}) 
(x^{T}x)^{d/2 - 1}u_{k}}
{(x^{T}x)^{d}}\\
&=& {\cal A}_{(k)} \tilde{u}_{k} \:-\: d \cdot (u_{k}^{T}{\cal A}_{(k)} \tilde{u}_{k}) u_{k}
\rule{0pt}{16pt}
\end{eqnarray*}
completing the proof of the theorem \end{proof}

It turns out that if the gradient of $\phi_{\cal A}^{(sym)}(x)$ is zero, then
the vector $x$ generally has the property that each subvector $u_{k}$ has the same norm.

\medskip

\begin{corollary}\label{symnab}
If $\nabla \phi_{\cal A}^{(sym)}(x) = 0$ and $x^Tx=1$, then either ${\cal A}_{(k)}\tilde{u}_{k} = 0$
for $k=1:d$ or
\[
\left[
\begin{array}{c}
{\cal A}_{(1)} \tilde{u}_{1} \\ \vdots \\ {\cal A}_{(d)} \tilde{u}_{d}
\end{array} \right]
\;=\;
d \cdot \phi_{\cal A}^{(sym)}(x)
\left[ \begin{array}{c}
u_{1} \\ \vdots \\
u_{d}
\end{array}
\right]
\]
and $\norm{u_1}_2 = \norm{u_2}_2 = \cdots = \norm{u_d}_2 = 1/\sqrt{d}$.
\end{corollary} 
\begin{proof}
Since $\nabla \phi_{\cal A}^{(sym)}(x) = 0$, we know from Theorem 4.5 that
\[
{\cal A}_{(k)} \tilde{u}_{k} \;=\; d \cdot (u_{k}^{T} {\cal A}_{(k)}\tilde{u}_{k}) u_{k}
\]
for $k=1:d$. Thus,
\[
u_{k}^{T}{\cal A}_{(k)} \tilde{u}_{k} \;=\; d \cdot (u_{k}^{T} {\cal A}_{(k)}\tilde{u}_{k}) (u_{k}^{T}u_{k}).
\]
From Lemma 3.1, if $u_{k}^{T}{\cal A}_{(k)} \tilde{u}_{k} = 0$ for some $k$, then it is
zero for all $k$. In this case we conclude from \eqref{phisymgrad} that ${\cal A}_{(k)} \tilde{u}_{k} = 0$
for $k=1:d$. Otherwise, $1 = d u_{k}^{T}u_{k}$, $k=1:d$. It follows that
$\norm{u_{1}}_{2}\; =\; \cdots\; =\; \norm{u_{d}}_{2} \;=\; 1/\sqrt{d}$.
\qquad \end{proof}

We are now ready for the main result that relates the eigenvalues and vectors of $\mathbf{sym}({\cal A})$ to the singular values and vectors of $\cal A$.

\begin{theorem}
If $\sigma $ is a nonzero singular value of $\inTens{A}{}{n_{1}}{n_{d}}$ with unit modal singular vectors
$u_{1},\ldots, u_{d}$, then 
\[
x_{\alpha} 
\;=\; \frac{1}{\sqrt{d}}
\left[ \begin{array}{c} 
 \alpha_{1}u_1 \\ 
  \alpha_{2}u_2\rule{0 pt}{10 pt} \\ 
\vdots \rule{0 pt}{10 pt}\\ 
\alpha{_d} u_d \rule{0 pt}{10 pt}
\end{array} \right] 
 \qquad \alpha = \left[ 1,\pm 1,\ldots,\pm 1\right]
\]
is an eigenvector for $\mathbf{sym}({\cal A})$ corresponding to eigenvalue
\[
\lambda_{\alpha}  \;=\; \alpha_{1}\alpha_{2}\cdots \alpha_{d}\: \frac{d!}{\sqrt{d^{d}}} \sigma.
\]
Note that $\alpha_1$ is set to +1 to resolve a uniqueness issue. See discussion after definition 4.4 and also equation (1.2) for the matrix case.
\end{theorem}

\medskip

\begin{proof}
We must show that $g = \nabla \phi_{\cal A}^{(sym)}(x_{\alpha}) = 0$.
If $\tau = \alpha_{1}\cdots \alpha_{d}$ then for $k=1:d$ we have from (4.6) that
\[
g_{k} \;= \frac{\tau}{\alpha_{k}d^{(d-1)/2}}{\cal A}_{k}\tilde{u}_{k}\:-\:
d \frac{\alpha_{k} \tau}{d^{(d+1)/2}} \left(u_{k}^{T}{\cal A}_{(k)} \tilde{u}_{k}\right)u_{k}
\]
But since $\sigma = u_{k}^{T}{\cal A}_{k}\tilde{u}_{k}$ and  ${\cal A}_{k}\tilde{u}_{k} = \sigma u_{k}$,
we have
\[
g_{k} \;=\; \alpha_{k} \tau \left( \frac{1}{d^{(d-1)/2}} \:-\: \frac{1}{d^{(d-1)/2}}\right) u_{k} \;=\; 0.
\]
Since 
$\lambda_{\alpha} \;=\; x_{\alpha}^{T} {\cal C}_{(1)}x_{\alpha}^{\otimes (d-1)}$
we have from Lemma 3.3 that
\begin{eqnarray*}
\lambda_{\alpha} &=&
\left(
\frac{1}{\sqrt{d}} \left[ \begin{array}{c} \alpha_{1}u_{1} \\ \vdots \\ \alpha_{d}u_{d} \end{array}\right]
\right)^{T}
\left(
\frac{(d-1)!}{d^{(d-1)/2}} 
\left[
\begin{array}{c}
(\tau/\alpha_{1}){\cal A}_{(1)} \tilde{u}_{1}\\
\vdots \\
(\tau/\alpha_{d}){\cal A}_{(d)} \tilde{u}_{d}
\end{array}
\right]
\right)\\
&=& \rule{0pt}{30pt}
\frac{1}{\sqrt{d}} \,\cdot \,
\frac{(d-1)!}{d^{(d-1)/2}}\,\cdot \,\tau \,\cdot \sum_{k=1}^{d} u_{k}^{T} {\cal A}_{(k)} \tilde{u}_{k}
\;=\;
\frac{(d-1)!}{\sqrt{d^{d}}} \,\cdot \, \tau \, \cdot \sum_{k=1}^{d} \sigma \;=\;
\frac{d!}{\sqrt{d^{d}}} \,\cdot \,\tau \,\cdot \, \sigma
\end{eqnarray*}
completing the proof of the theorem.
\end{proof}

\medskip

\noindent
Thus, for each singular value and vector for ${\cal A}$ we have $2^{d-1}$ eigenvalue/eigenvector
pairs for $\mathbf{sym}({\cal A})$.

\subsection{Connections to the Multilinear Transform}
Suppose $\inTens{F}{}{s_{1}}{s_{d}}$ and \inm{B_{k}}{s_{k}}{t_{k}} for $k=1:d$. The
tensor
$\inTens{T}{}{t_{1}}{t_{d}} $
defined by
\begin{equation} \label{mltransform}
{\cal T}(\mathbf{i}) = \sum_{\mathbf{j}=\mathbf{1}}^{\mathbf{s}} {\cal F}(\mathbf{j}) B_1(j_1,i_1)B_2(j_2,i_2)\cdots B_k( j_k,i_k).
\end{equation}
is the \emph{multilinear transform} \cite{lim} of tensor $\cal F$ by the matrices $B_{1},\ldots, B_{d}$ and is denoted by
\begin{equation}\label{mltransdot}
{\cal T} \:=\: {\cal F}\cdot(B_{1},B_2,\ldots,B_{d}).
\end{equation}
We also define
\begin{equation} \label{mltransdot2}
(B_1,B_2,\ldots,B_d) \cdot \mathcal{F} \;\; \equiv\;\;  \mathcal{F} \cdot (B_1^T,B_2^T,\ldots,B_d^T).
\end{equation}
Some of the key summations and vectors above can be expressed neatly through this transformation.
For example, if $\inTens{A}{}{n_{1}}{n_{d}}$ and \inv{u_{k}}{n_{k}} for $k=1:d$, then 
\[
{\cal A}\cdot (u_{1},\ldots,u_{d}) \;=\; 
\sum_{\mathbf{i}=\mathbf{1}}^{\mathbf{n}} {\cal A}(\mathbf{i})u_{1}(i_{1})\cdots u_{d}(i_{d}) = u_1^T \mathcal{A}_{(1)} \tilde{u}_1
\]
and 
\[
{\cal A}\cdot(u_1,\ldots,u_{k-1},I_{n_k},u_{k+1},\ldots,u_d) \;=\;
{\cal A}_{(k)} \tilde{u}_k.
\]

\section{Higher Order Power Methods}
We now briefly review various tensor power methods and consider them in light of the singular- and eigenvalue connection between $\cal A$ and $\mathbf{sym}(\mathcal{A})$.
%

%

\subsection{The HOPM}
The matrix power method method can be generalized to tensors by replacing the matrix-vector multipication with multilinear transforms. The \emph{Higher-Order Power Method} of \cite{lath3,lath4} for finding a singular value and associated singular vectors of general order-$d$ tensors proceeds in an alternating fashion to update each of the mode-$j$ singular vectors $u_j$.

\begin{algorithm}
  \caption{The higher-order power method (HOPM) \cite{lath3,lath4} }
  \label{alg:hopm}
    Given an order-$d$ tensor $\mathcal{A} \in \RR^{n_1\times \cdots \times n_d}$.
  \begin{algorithmic}[1]
    \Require ${u}_j^{(0)} \in \RR^{n_j}$ with $\|{u}_j^{(0)} \|_2 = 1$. Let
    $\sigma^{(0)} = (u_1^{(0)})^T\mathcal{A}_{(1)} \tilde{{u}}_1^{(0)}$.
    \For{$k=0,1,\dots$}
     \For{$j=0,1,\dots$}
      \State $\hat{{u}}_j^{(k+1)} \gets \mathcal{A}_{(j)} u_d^{(k+1)}\otimes \cdots \otimes u_{j+1}^{(k+1)} \otimes u_{j-1}^{(k)}\otimes \cdots \otimes u_1 ^{(k)}$
      \State $u_j^{(k+1)}  \gets \hat{{u}}_j^{(k+1)} /  \| \hat{{u}}_j^{(k+1)} \|_2$
     \EndFor
    \State $\sigma^{(k+1)} \gets (u_1^{(k+1)})^T \mathcal{A}_{(1)}\tilde{{u}}_1^{(k+1)}$
    \EndFor
  \end{algorithmic}
\end{algorithm}

%

\noindent
Different initial values for the $u_j$ vectors will in general result in convergence to different singular values. See Section 5.4 for a discussion on popular choices for higher-order power method initial values.

The HOPM can also be viewed as a way of finding the best rank-1 tensor approximation $\hat{\cal A}$ to $\cal A$  \cite{lath3}.
Specifically, a tensor $\inTens{\cal T}{}{n_1}{n_d}$ is said to be  \emph{rank-1} if for $k=1:d$ there exist vectors $\inv{t_i}{n_i}$ such that for all $\mathbf{i}=\mathbf{1},\ldots,\mathbf{n}$
\begin{equation}\label{r1tenseq1}
\mathcal{T}(\mathbf{i}) = t_1(i_1) t_2(i_2) \cdots t_d(i_d)
\end{equation}
and we then say that $\cal T$ is the \emph{tensor outer product} of the vectors $t_1,\ldots,t_d$, denoted by
\begin{equation}\label{r1tenseq2}
\mathcal{T} = t_1 \circ t_2 \circ \cdots \circ t_d.
\end{equation}
It can be shown  that the HOPM converges to a local minimum of the functional
$
f(\hat{\mathcal{A}}) \equiv \norm{\mathcal{A}-\hat{\mathcal{A}}}_F^2,
$
where $\hat{\cal A}= \sigma \,u_1 \circ\, \cdots \,\circ u_d$ is a rank-1 approximation to $\cal A$ and the Frobenius norm of a tensor $\cal T$ is defined as $\norm{\mathcal{T}}_F \equiv \sqrt{\sum_{\mathbf{i}=\mathbf{1}}^{{\mathbf{n}}} \mathcal{T}(\mathbf{i})^2}$. See \cite{regalia2}.


The HOPM can be applied to an order-$d$ symmetric $N\times \cdots \times N$ tensor, starting with a symmetric initial guess $u_1^{(0)} = u_2^{(0)} = \cdots = u_d^{(0)} \in \RR^N$. The solution found by the algorithm will be symmetric but intermediate results may break symmetry. Indeed, after one iteration the $u_j$ vectors will in general all be distinct, but $u_j^{(k)} \rightarrow u$ as $k \rightarrow \infty$ for some $u \in \RR^N$ \cite{lath3}.

\subsection{The S-HOPM}
Recently, \cite{regalia} investigated a modified version of the HOPM for symmetric tensors which was originally dismissed by \cite{lath3} as unreliable since in general it is not guaranteed to converge. This algorithm is called the \emph{Symmetric Higher Order Power Method} (S-HOPM) and converges for certain classes of symmetric tensors. For example, suppose $\cal C$ is a symmetric tensor of even order and that $M$ is a  square unfolding of $\cal C$. If $M$ is semidefinite then the S-HOPM converges \cite{regalia}.


\begin{algorithm}
  \caption{Symmetric higher-order power method (S-HOPM) \cite{lath3, regalia}}
  \label{alg:shopm}
    Given an order-$d$ symmetric tensor $\mathcal{C} \in \RR^{N\times \cdots \times N}$.
  \begin{algorithmic}[1]
    \Require ${x}^{(0)} \in \RR^N$ with $\|{x}^{(0)} \|_2 = 1$. Let
    $\lambda^{(0)} = (x^{(0)})^T\mathcal{C}_{(1)} ({x}^{(0)})^{\otimes (d-1)}$.
    \For{$k=0,1,\dots$}
    \State $\hat{{x}}^{(k+1)} \gets \mathcal{C}_{(1)} ({x}^{(k)})^{\otimes (d-1)}$
    \State ${x}^{(k+1)} \gets \hat{{x}}^{(k+1)} / \| \hat{{x}}^{(k+1)} \|_2$
    \State $\lambda^{(k+1)} \gets (x^{(k+1)})^T \mathcal{C}_{(1)}({x}^{(k+1)})^{\otimes (d-1)}$
    \EndFor
  \end{algorithmic}
\end{algorithm}


\noindent
This approach avoids the awkward situation, mentioned previously, of encountering non-symmetric intermediate values when using the HOPM on a symmetric tensor. 

Since $ \mathbf{sym}({\mathcal{A}})$ is symmetric for any tensor $\mathcal{A}$, the S-HOPM can be applied to $\mathcal{A}$ through its embedding.  By using facts previously established, we can reduce all operations on $ \mathbf{sym}({\mathcal{A}})$ to equivalent ones on $\mathcal{A}$. \myskip

\begin{algorithm}
  \caption{Symmetric higher-order power method on $\mathbf{sym}(\mathcal{A})$}
  \label{alg:symshopm}
    Given an order-$d$ tensor $\mathcal{A} \in \RR^{n_1\times \cdots \times n_d}$.
  \begin{algorithmic}[1]
    \Require ${u}_j^{(0)} \in \RR^{n_j}$ with $\|{u}_j^{(0)} \|_2 = 1$. Let
    $\sigma^{(0)} = (u_1^{(0)})^T\mathcal{A}_{(1)} \tilde{{u}}_1^{(0)}$.
    \For{$k=0,1,\dots$}
     \For{$j=0,1,\dots$}
      \State $\hat{{u}}_j^{(k+1)} \gets \mathcal{A}_{(j)} \tilde{{u}}_j^{(k)}$
      \State $u_j^{(k+1)}  \gets \hat{{u}}_j^{(k+1)} /  \| \hat{{u}}_j^{(k+1)} \|_2$
     \EndFor
    \State $\sigma^{(k+1)} \gets (u_1^{(k+1)})^T \mathcal{A}_{(1)}\tilde{{u}}_1^{(k+1)}$
    \EndFor
  \end{algorithmic}
\end{algorithm}


\noindent
This algorithm computes a singular value $\sigma$ for $\mathcal{A}$ and the mode-$j$ singular vectors $u_j$. The normalization used in Algorithm 3 is slightly different than a direct application of the S-HOPM on $\mathbf{sym}(\mathcal{A})$ would imply; the S-HOPM would set  $u_j^{(k+1)} =\hat{{u}}_j^{(k+1)} / \sqrt{ \lVert\hat{u}_1^{(k+1)} \rVert_2^2 +\cdots + \lVert\hat{u}_d^{(k+1)} \rVert_2^2}$. However, numerical experiments suggest that using $u_j^{(k+1)} =\hat{u}_j^{(k+1)} /  \lVert\hat{u}_j^{(k+1)} \rVert$ improves convergence. If $\cal A$ is itself symmetric, then Algorithm 3 reduces to the S-HOPM as all the $u_j$ will be equal, assuming $u_1^{(0)}=\cdots = u_d^{(0)}$. 

Note that Algorithm 3  is very similar to the regular HOPM except the most recently available information on $u_1,\ldots,u_{j-1}$ is not used when computing $u_j^{(k+1)}$ for $j>1$. The difference between the HOPM and Algorithm 3 is thus somewhat like the difference between the Jacobi and Gauss-Seidel iterative linear system solvers \cite{gvl}. 

%
%

Unlike the HOPM, Algorithm 3 does not always converge and since it can be shown that a square unfolding of $\mathbf{sym}(\mathcal{A})$ is indefinite unless all the entries in ${\cal A}$ are zero, the convergence criteria in \cite{regalia} do not apply.

%
%

\subsection{The SS-HOPM and $\mathbf{sym}(\cdot)$}Recently, Kolda and Mayo \cite{kolda4} developed a shifted version of the S-HOPM and proved that for a suitable choice of shift their algorithm will converge to an eigenpair $(\lambda,{x})$ for any symmetric tensor $\cal C$.

\begin{algorithm}
  \caption{Shifted symmetric higher-order power method (SS-HOPM) \cite{kolda4}}
  \label{alg:sshopm}
    Given an order-$d$ symmetric tensor $\mathcal{C} \in \RR^{N\times \cdots \times N}$.
  \begin{algorithmic}[1]
    \Require ${x}^{(0)} \in \RR^N$ with $\|{x}^{(0)} \|_2 = 1$. Let
    $\lambda^{(0)} = (x^{(0)})^T\mathcal{C}_{(1)} ({x}^{(0)})^{\otimes (d-1)}$.
    \For{$k=0,1,\dots$}
    \State $\hat{{x}}^{(k+1)} \gets \mathcal{C}_{(1)} ({x}^{(k)})^{\otimes (d-1)} + \alpha_\mathcal{C} x^{(k)}$
    \State ${x}^{(k+1)} \gets \hat{{x}}^{(k+1)} / \| \hat{{x}}^{(k+1)} \|_2$
    \State $\lambda^{(k+1)} \gets (x^{(k+1)})^T \mathcal{C}_{(1)}({x}^{(k+1)})^{\otimes (d-1)}$
    \EndFor
  \end{algorithmic}
\end{algorithm}

%

\noindent
If the shift $\alpha_\mathcal{C}$ satisfies $\lvert \alpha_\mathcal{C} \rvert > (d-1) \sum_{\mathbf{i} = \mathbf{1}}^{\mathbf{N}} \lvert \mathcal{C}(\mathbf{i}) \rvert$ then the SS-HOPM will converge to an eigenpair \cite{kolda4}.

When ${\cal C} = \mathbf{sym}({\cal A})$ the algorithm can be simplified and expressed in terms of operations on $\cal A$. 

\begin{algorithm}
  \caption{Shifted symmetric higher-order power method on $\mathbf{sym}(\mathcal{A})$}
  \label{alg:symsshopm}
    Given an order-$d$ tensor $\mathcal{A} \in \RR^{n_1\times \cdots \times n_d}$.
  \begin{algorithmic}[1]
    \Require ${u}_j^{(0)} \in \RR^{n_j}$ with $\|{u}_j^{(0)} \|_2 = 1/\sqrt{d}$. Let
    $\sigma^{(0)} = (u_1^{(0)})^T\mathcal{A}_{(1)} \tilde{{u}}_1^{(0)}$.
    \For{$k=0,1,\dots$}
     \For{$j=0,1,\dots$}
      \State $\hat{{u}}_j^{(k+1)} \gets \mathcal{A}_{(j)} \tilde{{u}}_j^{(k)} + d\, \alpha_\mathcal{A} u_j^{(k)}$
     \EndFor
     \For{$j=0,1,\dots$}
      \State $u_j^{(k+1)}  \gets \hat{{u}}_j^{(k+1)} / \sqrt{ \| \hat{{u}}_1^{(k+1)} \|_2^2 + \dots + \| \hat{u}_d^{(k+1)} \|_2^2 }$
     \EndFor
    \State $\sigma^{(k+1)} \gets (u_1^{(k+1)})^T \mathcal{A}_{(1)}\tilde{{u}}_1^{(k+1)}$
    \EndFor
  \end{algorithmic}
\end{algorithm}
%
%
%

\noindent
Using Theorem 4.7, a simple normalization of the values returned by this algorithm gives a singular value and associated unit norm singular vectors of $\cal A$.

The shift $\alpha_\mathcal{A}$ must satisfy  $\lvert \alpha_\mathcal{A} \rvert > (d-1) \sum_{\mathbf{i} = \mathbf{1}}^{\mathbf{n}}\lvert \mathcal{A}(\mathbf{i})\rvert$ to guarantee convergence, although a smaller shift might be sufficient for any particular tensor $\cal A$. 
\medskip

\paragraph{Example}Let $\cal A$ be the $2\times 2 \times 2 \times 2$ tensor given by the unfolding
\[
\mathcal{A}_{(1)} =
\small{ \left[ \! \begin{array}{rrrrrrrr}
    1.1650 &   0.2641 &  -0.6965 &   1.2460 &   0.0751 &  -1.4462 &   0.0591&    0.5774\\
    0.6268 &   0.8717 &   1.6961&   -0.6390 &   0.3516 &  -0.7012 &   1.7971&   -0.3600
    \end{array}\! \right]. }
\]
We ran 100 trials of the HOPM, Algorithm 3 and Algorithm 5 using different random starting points $u_i^{(0)}$ chosen from a uniform distribution on $[-1, 1]^{n_i}$ and suitably normalized for each algorithm. The algorithms are considered to have converged when $\lvert \sigma^{(k+1)} - \sigma^{(k)} \rvert < 10^{-16}$. For this example, all three algorithms converged for every starting point.

The HOPM found the singular values $2.7248$ and $1.7960$. Algorithm 3 converged to $\sigma = \pm 2.7248$. Algorithm 5 with a positive shift $\alpha_{\cal A}$ found $2.7248$ and $1.7960$ and using a negative shift produced the values $-2.7248$ and $-1.7960$.

For this tensor $\cal A$ the theory suggests a shift $\alpha_{\cal A}$ greater than 37.72 in absolute value to guarantee convergence. However, using $\alpha_{\cal A}$ as small as 1 will still lead to convergence and does so in many fewer iterations, sometimes by as much as a factor of 30 when compared to the suggested shift. Setting $\alpha_{\cal A}$ to zero caused the algorithm to fail to converge for all chosen starting points.

\subsection{Initialization}A standard way to initialize higher-order power methods is to use a truncated form of the Higher-Order Singular Value Decomposition (HOSVD) of \cite{lath},
\begin{equation}\label{hosvd1}
\mathcal{A} =  ( U_1,U_2,\ldots, U_d) \cdot \mathcal{S}.
\end{equation}
where $\inTens{\cal S}{}{n_1}{n_d}$ is the \emph{core tensor}, the $U_i  \in \RR^{n_i \times n_i}$ are orthogonal and related to the modal unfoldings of $\cal A$ through the matrix SVD equations $\mathcal{A}_{(i)} \;=\; U_i \Sigma_i V_i^T$. 

To initialize the HOPM, for example, the values $u_j^{(0)} = U_j(:,1)$ have been shown \cite{lath3} to often lie close to the best rank-1 approximation to $\cal A$.

If desired, it is possible to create the HOSVD of ${\cal C} = \mathbf{sym}({\cal A})$ from the HOSVD of $\cal A$. For example, if $\mathcal{C} = (U_\mathcal{C},\ldots,U_\mathcal{C}) \cdot \mathcal{S}_\mathcal{C}$ is the HOSVD of $\cal C$ then it can be shown that $U_\mathcal{C}$ is a column permutation of the block-diagonal matrix $\mathbf{diag}(U_1,U_2,\ldots,U_d)$.

There are many other ways to initialize tensor power methods. In \cite{regalia2} Regalia and Kofidis derive a procedure for symmetric tensors that can outperform the HOSVD-based approach. 

Another possibility is to compute a tensor generalization of the QR decomposition with partial pivoting, of the form $\mathcal{A} = (Q_1,\ldots,Q_d) \cdot \mathcal{R}$ where $\mathcal{A}_{(k)} = Q_k R_k \Pi_k$ are the pivoted QR decompositions of the unfolding ${\cal A}_{(k)}$. It can be shown that this ``HOQRD'' decomposition retains some of the approximation properties of the truncated matrix pivoted QR decomposition and can thus give a reasonable initial guess for a tensor power method. As for the HOSVD, the HOQRD of $\mathbf{sym}(\mathcal{A})$ can be constructed from the HOQRD of $\cal A$.

\section{Tensor Rank and the {\bf sym} Operation}

There are several definitions of tensor rank, each of which represents some reasonable generalization of matrix rank. For an  excellent review see \cite{lim}. 
In this brief section we relate the multilinear rank and the outer product rank of 
$\mathbf{sym}({\cal A})$ to 
 the multilinear rank and the outer product rank of $\cal A$.

\subsection{Multilinear Rank}

The \emph{multilinear rank} of $\inTens{A}{}{n_{1}}{n_{d}}$ is the $d$-tuple $\mlrank{\mathcal{A}} = (r_1(\mathcal{A}),r_2(\mathcal{A}),\ldots,r_d(\mathcal{A}))$ where $r_i(\mathcal{A}) = \rank (\mathcal{A}_{(i)})$. Note that if the tensor $\inTens{C}{}{N}{N}$ is symmetric, and $R = \mbox{rank}({\cal C}_{(1)})$, then 
 \begin{equation}
\mathrm{rank}_{\boxplus}(\mathcal{C}) \;=\; (R,R,\ldots,R)
\end{equation}
because $\mathcal{C}_{(1)} = \mathcal{C}_{(2)} =\cdots = \mathcal{C}_{(d)}$. If 
$\mathcal{C} = \mathbf{sym}(\mathcal{A})$, then it is possible  to connect
$\mathrm{rank}_{\boxplus}(\mathcal{C})$ to $\mathrm{rank}_{\boxplus}(\mathcal{A})$.

\medskip

\begin{theorem}
If $\inTens{A}{}{n_{1}}{n_{d}}$
and $\mathrm{rank}_{\boxplus}(\mathcal{A}) = (r_1,\ldots, r_d)$,
then 
\[
\mathrm{rank}_{\boxplus}({\mathbf{sym}({\mathcal{A}})}) = (R,\ldots,R)
\]
where $R = r_{1} + \cdots + r_{d}$.
\end{theorem}
\begin{proof} Suppose $\mathcal{C} = \mathbf{sym}(\mathcal{A})$ and $\mathcal{C}_{\mathbf{i}}$
is $\mathcal{C}$'s $\mathbf{i}$th block, $\mathbf{1} \leq \mathbf{i} \leq \mathbf{d}$.
 Let $\mathcal{C}^{(k)}$ be a 
$1\times d \times \cdots \times d$ block tensor defined by
\[
\mathcal{C}_{\mathbf{i}}^{(k)} \;=\; \mathcal{C}_{\mathbf{i}}
\] 
where $i_{1} = k$ and $1 \leq i_{j} \leq d$ for $j=2:d$. Note that if $\mathbf{i}(2:d)$
is a permutation of $[1:k-1 \;\; k+1:d]$, then 
\[
\mathcal{C}_{\mathbf{i}}^{(k)} \;=\; \mathcal{A}^{<[\,k \:\mathbf{i}(2:d)\,]>}.
\]
It follows that 
\begin{equation}
\mbox{range}(\mathcal{C}_{(1)}^{(k)}) = \mbox{range}(\mathcal{A}_{(k)}) \qquad k=1:d
\end{equation}
If 
\[
\mathcal{C}_{(1)} \;=\; \left[ \begin{array}{c} C_{1} \\ \vdots \\ C_{d} \end{array}
\right]
\!\!
\begin{array}{l} 
\left. \right\} n_{1} \\ \vdots \\
\left. \right\} n_{d}
\end{array}
\]
is a block row partitioning of $\mathcal{C}_{(1)}$, then $C_{k}$ is a column permutation
of $\mathcal{C}_{(1)}^{(k)}$ and so using (6.2) we have
\begin{equation}
\mbox{rank}(C_{k}) \;=\; \mbox{rank}(\mathcal{C}_{(1)}^{(k)})\;=\; r_{k}.
\end{equation}
If
\[
v \;=\;
\left[ \begin{array}{c} v_{1} \\ \vdots \\ v_{d} \end{array}
\right]
\!\!
\begin{array}{l} 
\left. \right\} n_{1} \\ \vdots \\
\left. \right\} n_{d}
\end{array}
\]
is a column of $\mathcal{C}_{(1)}$ then it is a 
a mode-1 fiber of $\mathcal{C}$ and thus
can ``pass through'' at most one $\mathcal{C}$-block having
an index that is a permutation of $1:d$. This means that at most one of $v$'s subvectors is
zero. It follows from (6.3)
that
\[
\mbox{rank}(\mathcal{C}_{(1)}) \;=\; \sum_{k=1}^{d} \mbox{rank}(C_{k}) \;=\; \sum_{k=1}^{d}r_{k}
\]
completing the proof of the theorem.
\end{proof}

\subsection{Outer Product Rank}
The \emph{outer product rank} of $\inTens{\cal A}{}{n_1}{n_d}$ is the minimum number of rank-1 tensors, as defined in \eqref{r1tenseq1} and \eqref{r1tenseq2}, that are needed to represent it as a sum
$$ \oprank{\mathcal{A}} \equiv \min \left\{ r \; : \; \mathcal{A} = \sum_{i=1}^r u_1^{(i)} \circ u_2^{(i)} \circ \cdots \circ u_d^{(i)}, \quad u_j^{(i)}\in \RR^{n_j} \right\}. $$
For matrices $\rank (\mathbf{sym}({A})) = 2\, \rank(A)$. 
Indeed, If $A = \sum_{i=1}^r  \sigma_i u_i v_i^T$  is the SVD of $A$,then
 $$\mathbf{sym}({A}) 
 \;=\;
  \sum_{i=1}^{r}  \sigma_i \left(    
   \left[ \begin{array}{c} 0 \\ v_i \end{array}\right]
    \left[ \begin{array}{cc} u_i^T & 0 \end{array} \right] +
     \left[ \begin{array}{c} u_i \\ 0 \end{array}\right]
      \left[ \begin{array}{cc} 0 & v_i^T \end{array}\right]    \right).$$
Motivated by this expansion we make a definition.
\begin{definition}
If $\inTens{T}{}{n_{1}}{n_{d}}$ is the rank-1 tensor $\mathcal{T} = t_1 \circ \cdots \circ t_d$ 
and $N = n_{1}+\cdots n_{d}$, then $\inTens{S}{}{N}{N}$ is the rank-1 tensor
\[
\mathcal{S} \;=\; \pi(\mathcal{T})  \;=\; s_1 \circ \cdots \circ s_d
\]
where

\[
s_{k} \;=\; \left[ \begin{array}{c} 0 \\ t_{k} \rule{0pt}{12pt}
\\ 0 \rule{0pt}{12pt} \end{array} \right]
\!\! \begin{array}{l}
\left. \right\} n_{1} + \cdots + n_{k-1} \\
\left. \right\} n_{k} \rule{0pt}{12pt} \\
\left. \right\} n_{k+1} + \cdots + n_{d} \rule{0pt}{12pt}
\end{array}
\]
\end{definition}

\noindent
With this construction, we can produce an outer product expansion of
$\mathbf{sym}(\mathcal{A})$ given an outer product expansion of $\mathcal{A}$.
\begin{theorem}
If \inTens{A}{}{n_{1}}{n_{d}} and
\[
\mathcal{A}\; =\; \sum_{i=1}^r u_{i}^{(1)} \circ u_i^{(2)} \circ \cdots \circ u_{i}^{(d)}
\]
where \inv{u_{1}^{(k)},\ldots,u_{r}^{(k)}}{n_{k}}, then
\begin{equation}
\mathbf{sym}({\mathcal{A}}) \;=\; \sum_{\mathbf{p}\in S_{d}} \sum_{i=1}^r \pi( u_{i}^{(1)} \circ u_i^{(2)} \circ \cdots \circ u_{i}^{(d)} )^{<\mathbf{p}>}
\end{equation}
where $S_{d}$ is the set of all permutations of $1:d$.
\end{theorem}

\begin{proof}
Let $\mathcal{C}$ be the sum on the right side of (6.4)
and note that
\[
\mathcal{C} \;=\;
\sum_{\mathbf{p}\in S_{d}} \sum_{i=1}^r \pi( u_{i}^{(p_{1})} \circ u_i^{(p_{2})} \circ \cdots \circ u_{i}^{(p_{d})} ).
\]
We must show that the 
$\mathbf{q}$th block of $\mathbf{sym}(\mathcal{A})$ equals the
$\mathbf{q}$th block of $\mathcal{C}_{\mathbf{q}}$. If $\mathbf{q}$ is not a permutation
of $1:d$, then these blocks are both zero. Otherwise
\[
 \mathcal{C}_{\mathbf{q}} \;=\; \sum_{i=1}^r  u_{i}^{(q_{1})} \circ u_i^{(q_{2})} \circ \cdots \circ u_{i}^{(q_{d})} 
\;=\; 
\left(\sum_{i=1}^r  u_{i}^{(1)} \circ u_i^{(2)} \circ \cdots \circ u_{i}^{(d)} \right)^{<\mathbf{q}>}
\;=\; 
\mathcal{A}^{<\mathbf{q}>}.
\]
completing the proof of the theorem.
\end{proof}

\medskip

\noindent
Since the double summation in (6.4) involves $rd!$ terms, it follows that
\begin{equation}
 \mathrm{rank}_\otimes ({\mathbf{sym}({\mathcal{A}}})) \; \leq \;d! \cdot \mbox{rank}_{\otimes}(\mathcal{A})
\end{equation}
We conjecture that equality prevails.
This is somewhat reminiscent of the direct sum conjecture \cite{directsum}, i.e.~that $\oprank{\mathcal{A}\oplus \mathcal{B}} = \oprank{\mathcal{A}} + \oprank{\mathcal{B}}$. Intuitively, $\mathbf{sym}({\cal A})$ contains $d!$ distinct copies of $\cal A$ in nonoverlapping index regions so if the matrix case were to generalize, any expansion of $\mathbf{sym}({\cal A})$ into a sum of $\leq d!r$ rank-1 terms could be reduced to (6.4) without adding terms, thus having exactly $d!r$ terms. We have so far been unable to prove this.
Note that it can be shown that 
\[
d \cdot \mbox{rank}_{\otimes}(\mathcal{A}) \; \leq \; \mathrm{rank}_\otimes ({\mathbf{sym}({\mathcal{A}}}))
\]
using Lemma 3.5 in \cite{lim}.


\section{Conclusions}

The symmetrization $\mathbf{sym}(\mathcal{A})$ can be used to connect algorithms for symmetric tensors and ones for general tensors. In this paper we have shown how algorithms such as the S-HOPM and SS-HOPM give rise to non-symmetric algorithms through the symmetrization in a way that preserves many convergence properties. In particular, the non-symmetric version of the SS-HOPM we derive is guaranteed to converge for an appropriately chosen shift $\alpha_{\cal A}$. Are there other tensor methods where the symmetrization could be used to spot new connections or derive useful algorithms? 

The rank properties of the symmetrization in some ways mirror the matrix case, but fundamental questions regarding the outer product rank of $\mathbf{sym}(\mathcal{A})$ remain open. Resolution
of these questions may help bridge the conceptual gap that exists between matrix rank
 and tensor rank.

\section{Acknowledgements}

The authors  thank the referees for suggestions that led
to an improved presentation.

\end{document}